\documentclass[12pt]{article}
\usepackage{amsmath}
\usepackage{graphicx}
\usepackage{enumerate}
\usepackage{natbib}
\usepackage{url} 
\usepackage{array}
\usepackage{mathrsfs}
\usepackage{amssymb,bbm}
\usepackage{amsfonts}
\usepackage{amscd}
\usepackage{color}
\usepackage{latexsym,bm}
\usepackage{epstopdf}
\usepackage{rotating}
\usepackage{multirow,bigstrut,booktabs}
\usepackage[colorlinks,citecolor=blue,urlcolor=blue]{hyperref}

\newcommand{\blind}{1}

\newtheorem{theorem}{Theorem}

\newtheorem{assumption}{Assumption}

\newtheorem{proposition}{Proposition}

\newcommand{\Var}{\mathrm{Var}}

\newcommand{\E}{\mathbb{E}}

\begin{document}

\def\spacingset#1{\renewcommand{\baselinestretch}%
{#1}\small\normalsize} \spacingset{1}


\if1\blind
{
  \title{\bf Bubble Modeling and Tagging: A Stochastic Nonlinear Autoregression Approach}
  \author{
  Xuanling Yang\textsuperscript{1}, 
  Dong Li\textsuperscript{2}, 
  and Ting Zhang\textsuperscript{3}\thanks{Zhang gratefully acknowledges the support of the U.S. NSF DMS-2412661.} \vspace{.2cm} \\
  \small
  \textsuperscript{1} College of Mathematics and Physics,
  Beijing University of Chemical Technology\\
  \small
  \textsuperscript{2} Department of Statistics and Data Science,
  Tsinghua University\\
  \small
  \textsuperscript{3} Department of Statistics, University of Georgia}
  \maketitle
} \fi

\if0\blind
{
  \bigskip
  \bigskip
  \bigskip
  \begin{center}
    {\LARGE\bf Bubble Modeling and Tagging: A Stochastic Nonlinear Autoregression Approache}
\end{center}
  \medskip
} \fi

\bigskip
\begin{abstract}
Economic and financial time series can feature locally explosive behavior when a bubble
is formed. The economic or financial bubble, especially its dynamics, is an intriguing topic that has been
attracting longstanding attention. To illustrate the dynamics of the local explosion itself, the paper
presents a novel, simple, yet useful time series model, called the stochastic nonlinear autoregressive model, which is
always strictly stationary and geometrically ergodic and can create long swings or persistence observed
in many macroeconomic variables. When a nonlinear autoregressive coefficient is outside of a certain range, 
the model has periodically explosive
behaviors and can then be used to portray the bubble dynamics. Further, the quasi-maximum likelihood
estimation (QMLE) of our model is considered, and its strong consistency and asymptotic normality are
established under minimal assumptions on innovation. A new model diagnostic checking statistic is
developed for model fitting adequacy. In addition, two methods for bubble tagging are proposed,
one from the residual perspective and the other from the null-state perspective.
Monte Carlo simulation studies are conducted
to assess the performances of the QMLE and the two bubble tagging methods in finite samples.
Finally, the usefulness of the model is illustrated by an empirical
application to the monthly Hang Seng Index.
\end{abstract}

\noindent%
\vfill

\newpage
\spacingset{1.9} 
\section{Introduction}\label{sec:intro}
Financial speculative bubbles have been attracting longstanding attention of economists and financial practitioners as an economic crisis often originates along with a burst of a bubble.
In reality, however, economic or financial bubbles cannot be avoided.
The presence of bubbles is partially evidenced by that many economic or financial time series possess locally explosive behavior and a subsequent burst, with such a phenomenon appearing periodically.
Studying the dynamics of bubble thus becomes important and intriguing.

One classical definition of the bubble is the deviation of the market price from its fundamental value (a sum of discounted future dividends) in rational expectation price models.
An important model of the rational bubble is initiated by
\cite{Blanchard1982Bubbles}, where the bubble process is captured via a simple stochastic autoregression (AR) with a fixed explosive rate and an absorbing state zero.
Their model is then extended by \cite{Evans1991Pitfalls} via adopting a stochastic rate of explosion.
Primarily, the bubble is regarded as an explosive nonstationary process,
which motivates to test its presence via unit root and cointegration tests \citep{Diba1988a,Diba1988b}.
Recently, this idea is further developed by \cite{Phillips2011Dating},
\cite{Phillips2011Explosive, Phillips2015Testing, Phillips2015},
 \cite{Harvey2019, Harvey2020}, \cite{Taoyu}, \cite{Kurozumi2023Time},
\cite{Esteve2023} and references therein. On the other hand, \cite{Evans1991Pitfalls} also notes that periodical collapse of bubbles makes the bubble paths look more like a stationary process.
Within a stationary framework, \cite{Gourieroux2017Local} find that noncausal AR(1) models can characterize multiple local explosions in time series. Then this noncausal approach
to bubble modelling has been extended to high-order
mixed causal-noncausal\footnote{The definition of `causal' in time series can be referred
to \cite{BrockwellDavis1991}. Early, \cite{Lanne2011Noncausal} study statistical inference on noncausal AR models with an application to U.S. inflation dynamics.} time series models, see, for example, \cite{Gouriroux2016FilteringPA}, \cite{Fries2019Mixed},
\cite{Cavaliere2018BootstrappingNA},
 \cite{Davis2020NoncausalVA}, and \cite{Fries2021ConditionalMO}.
However, one shortcoming of the noncausal approach invites computational challenge
and many resampling methods are needed.
To bypass this shortcoming, \cite{Blasques2022A} propose a new observation driven model with time-varying parameters
and study its probabilistic properties and statistical inference.
Nevertheless, their estimation heavily depends on the choice of the survival function and the asymptotics
can be obtained only for a part of parameters.
Motivated by all above facts, we here present a new simple time series model to describe the dynamics of bubbles.

In this paper, a first-order stochastic nonlinear autoregressive (SNAR) model $\{y_t\}$ is defined as
\begin{flalign}\label{eq.model}
y_t=s_t\phi_0 |y_{t-1}|+\varepsilon_t, \quad t\in\mathbb{Z}:=\{0, \pm1,\pm2, ...\},
\end{flalign}
where $\phi_0\in \mathbb{R}$, $\{\varepsilon_t: \, t\in\mathbb{Z}\}$ is a sequence of independent
and identically distributed (i.i.d.) random variables
on some basic probability space $(\Omega, \mathcal{F}, \mathbb{P})$,
and independent of i.i.d. binary variables $\{s_t: \,t\in\mathbb{Z}\}$ with
$\mathbb{P}(s_t=1)=p_0=1-\mathbb{P}(s_t=0)$, $p_0\in[0, 1)$.

\begin{figure}[!htbp]
\begin{center}
\includegraphics[width=15cm]{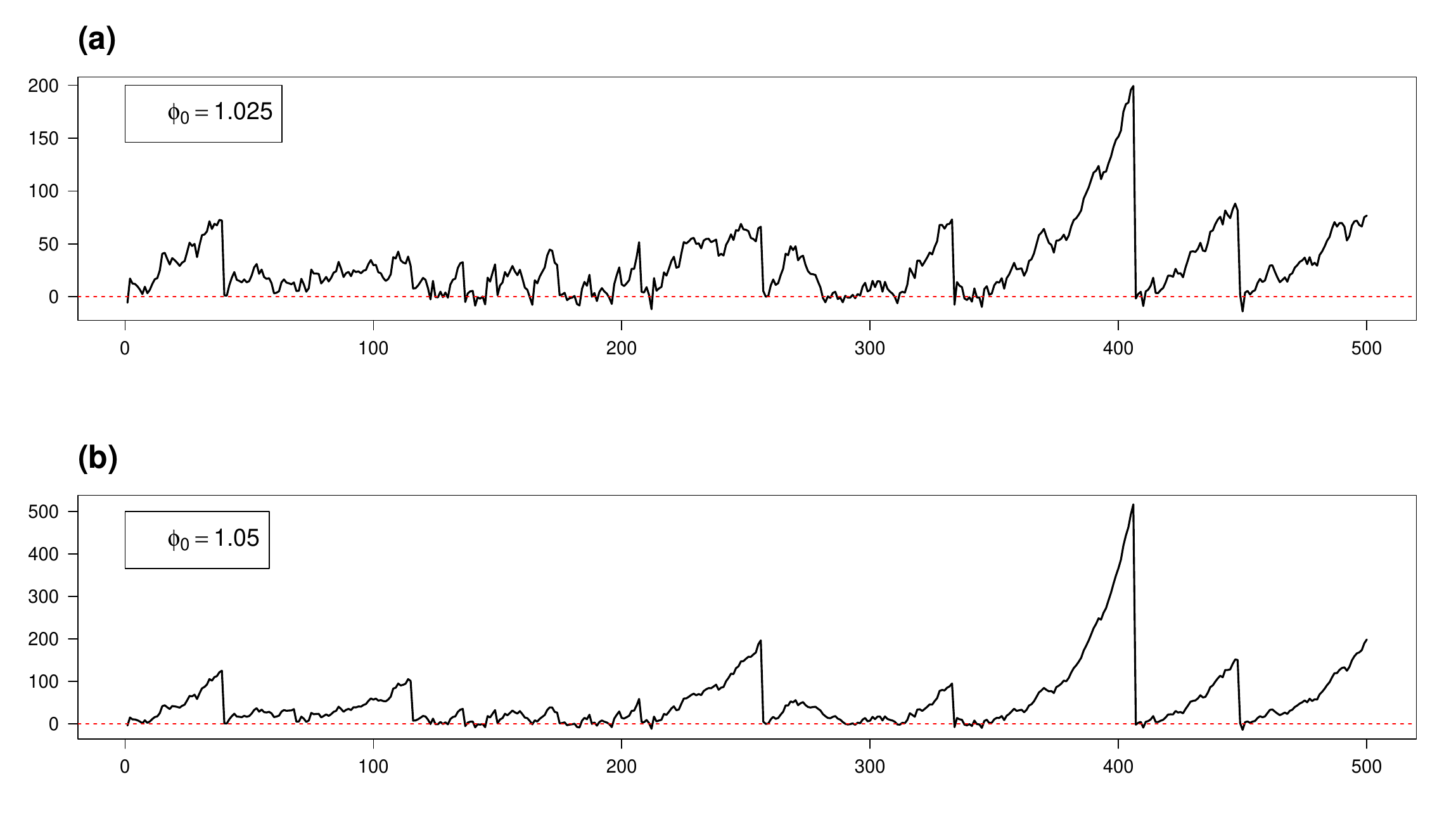}
\end{center}
\caption{Simulated paths of model \eqref{eq.model} with $\varepsilon_t\sim\mathcal N(0,6^2)$,  $p_0=0.977$, and (a) $\phi_0=1.025$ and (b) $\phi_0=1.05$. \label{fig.path}}
\end{figure}

Clearly, when $\phi_0>1$, $y_t$ is explosive in the periods where $s_t=1$
and creates an excursion which stops once $s_t=0$.
Fig.~\ref{fig.path} illustrates two simulated paths of the SNAR model (\ref{eq.model})
with $\varepsilon_t\stackrel{i.i.d}{\sim}\mathcal N(0,6^2)$.
We can observe periodically local explosions followed by bursts. 
And the 'shape' of a bubble before its burst is similar to a quadratic curve, which conforms to some practitioners' pointview that the accumulation of a bubble before bursting resembles a parabola.
Further, when $\mathbb{P}(s_t=1)=1$, the SNAR model (\ref{eq.model}) reduces
to an absolute AR model or a special threshold AR model with threshold parameter zero, which is studied in \cite{Tong1990}, \cite{LiTong2020} and references therein.

Major contributions of our paper are as follows.

First, we introduce a simple yet useful time series model, a SNAR model, for modeling the dynamics of bubble process.
We then prove that the model is always strictly stationary and geometrically ergodic under minimal assumptions on innovation and the probability $p_0$. Within a causal and stationary framework, when the parameter $\phi_0>1$,
our SNAR model still displays periodically local explosions and collapses. It can create long swings or persistence and can then be used to portray the bubble dynamics.
More importantly, our model is always causal in the classical sense of time series.
Compared with the noncausal bubble models in the literature,
our model facilitates the understanding of the dynamics of bubble,
and is quite simpler and more convenient in application. Additionally,
our model avoids computational burden of the noncausal approach and  keeps away from the choice of the survival function in the time-varying parameter model in \cite{Blasques2022A}.

It is worth mentioning that a related model to ours is a stochastic AR initiated by \cite{Blanchard1982Bubbles},\footnote{It is a simple case of a
random coefficient AR model in \cite{Nicholls1982RandomCA} if we regard  $s_t\phi_0$ as a coefficient.
Random coefficient AR models have received great attention and been widely used in econometrics, finance,
 engineering, among others, due to their flexibility, parsimonious representation and analytical tractability.} which is defined as
$y_t=s_t\phi_0 y_{t-1}+\varepsilon_t$, $t\in\mathbb{Z}$, where $\{s_t\}$ and $\{\varepsilon_t\}$ are defined in
(\ref{eq.model}). It can also create long swings or persistence \citep{Johansen2013}. Unfortunately,
it can usually generate many negative local explosions even if $\phi_0>1$ and
$\varepsilon_t\stackrel{i.i.d}{\sim}\mathcal N(0,1)$. This does not meet our common sense on financial bubble.

Second, we consider the quasi-maximum likelihood estimation (QMLE) of the SNAR model and establish
its strong consistency and asymptotic normality under minimal assumptions on innovation
and probability parameter $p_0$, regardless of infinite variance or heavy-tailedness of the model.

Third, we develop a new model diagnostic checking statistic since the classical portmanteau test is
invalid for our model.\footnote{The reason is that the residuals cannot be obtained.}

Fourth, we consider two methods for tagging the bubbles, one from the residual point of view and the other from the null-state perspective. The problem of bubble detection has been studied in the literature; see for example \cite{Phillips2011Dating}, \citet{Phillips2015Testing}, \cite{Phillips2015}, \cite{Blasques2022A}, \cite{Kurozumi}, and references therein. Existing results in this direction, however, have been mainly developed by viewing the bubble as a separate process occurring on an unknown but deterministic time interval within the observation period. The current paper, on the other hand, aims to incorporate the bubble mechanism into the data generating process to provide a stationary statistical model that can capture and interpret bubbles, including their formations and collapses. As a result, unlike existing results that typically assume the bubbles to persist for an adequate duration to achieve their consistent detection, the bubbles in the current model can be transient and thus the problem of bubble tagging can be more challenging in the current setting. For this, we consider two approaches, where the first one utilizes the nonlinear autoregressive residual from the proposed model and the second one is constructed from a hypothesis testing point of view. For both methods, we provide theoretical quantification on the finite-sample probability of correct tagging under reasonably mild conditions. Monte Carlo simulation results are provided to assess the finite-sample performance of the proposed QMLE and bubble tagging methods.

The remainder of the paper is organized as follows. Section \ref{stable} investigates strict stationarity and
geometric ergodicity of model (\ref{eq.model}).
Section \ref{section.qmle} considers the QMLE with its
asymptotics. Section \ref{section.checking} studies model diagnostic checking.
Section \ref{section.label} considers the problem of bubble tagging, where two approaches are considered with their finite-sample probability bounds studied.
Section \ref{section.simulation} carries out Monte Carlo simulation studies to assess the finite-sample performances of the QMLE and the two bubble tagging methods.
Section \ref{section.example} gives an empirical application to illustrate the usefulness of the model.
Section \ref{conclusion} concludes.  All technical proofs are relegated to the Supplementary Material.

\section{Probabilistic Properties of the SNAR Model}\label{stable}
The aim of this section is to prove the strict stationarity and geometric ergodicity of model
(\ref{eq.model}) under a very mild condition. We will prove the following result, using the approach
developed by \cite{MeynTweedie} for establishing the geometric ergodicity of Markov chains.
This result is important and is a theoretical foundation of statistical inference for model (\ref{eq.model}).
\begin{theorem}\label{ssge}
Suppose that $\mathrm{(i)}$ $\{\varepsilon_t\}$ is i.i.d. and independent of i.i.d. binary variables $\{s_t\}$ with $0\leq p_0<1$, and $\mathrm{(ii)}$ $\varepsilon_1$ has a positive density on $\mathbb{R}$ with $\mathbb{E}(\log^+|\varepsilon_1|)<\infty$.
Then there exists a strictly stationary, nonanticipative solution to $\{y_t\}$ in model $(\ref{eq.model})$
and the solution is unique and geometrically ergodic.
\end{theorem}

\begin{figure}[!htbp]
  \centering
  \includegraphics[width=12cm]{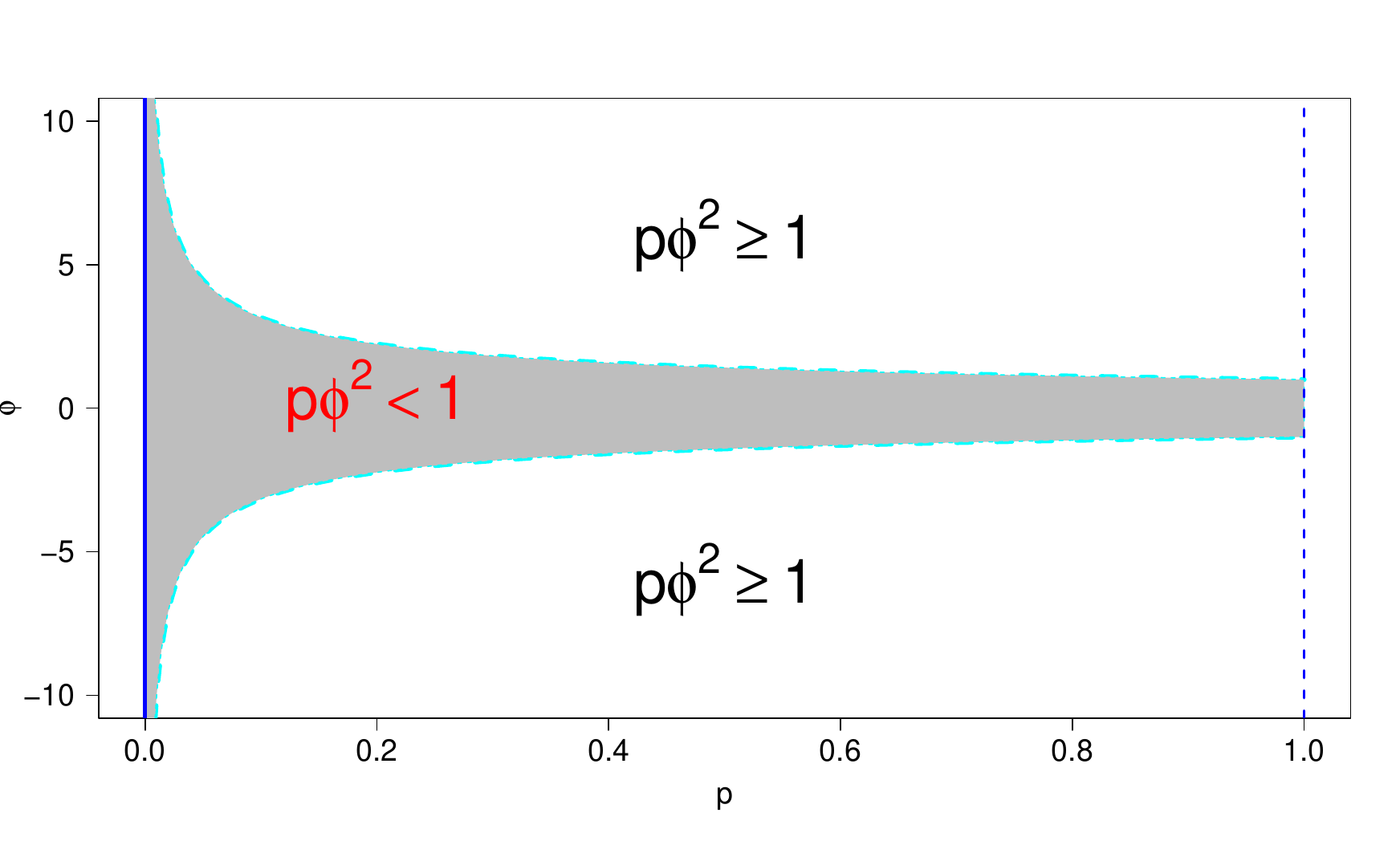}
  \caption{The strict stationarity region $\{(p, \phi): \,\phi\in\mathbb{R}, \,p\phi^2<1, \,0\leq p<1\}$ of  $y_t$ with finite second moment.}\label{fig.region}
\end{figure}

Next, consider the existence conditions on moments of $y_t$. Clearly, when $\mathbb{E}\varepsilon_t=0$, $\mathbb{E}\varepsilon_t^2<\infty$,
and both $\varepsilon_t$ and $s_t$ are independent,  then
$\mathbb{E}y_t^2=\mathbb{E}\varepsilon_t^2/(1-p_0\phi_0^2)<\infty$ if $p_0\phi_0^2<1$. Fig.~\ref{fig.region} plots the strict stationarity region of $y_t$ with $\mathbb{E}y_t^2<\infty$.
Further, if $\E \varepsilon_t^3=0$, then a simple algebraic calculation gives the kurtosis of $y_t$:
\begin{flalign*}
\mathrm{kurtosis}(y_t)=\frac{\left\{6p_0\phi_0^2+\mathrm{kurtosis(\varepsilon_t)}
(1-p_0\phi_0^2)\right\}(1-p_0\phi_0^2)}{1-p_0\phi_0^4},\quad\mbox{if}\quad p_0\phi_0^4<1.
\end{flalign*}
In particular, when $\varepsilon_t\sim \mathcal{N}(0, 1)$, then
\begin{flalign*}
\mathrm{kurtosis}(y_t)=\frac{3(1-p_0^2\phi_0^4)}{1-p_0\phi_0^4}>3,\quad\mbox{if}\quad 0<p_0\phi_0^4<1,
\end{flalign*}
which implies that $\{y_t\}$ must be heavy-tailed.

\section{Quasi-Maximum Likelihood Estimation}\label{section.qmle}
Let $\theta_0=(\phi_0, p_0, \sigma_0^2)'$ be the true parameter with
$\sigma^2_0=\mathbb{E}\varepsilon_t^2$. Denote by $\theta=(\phi, p, \sigma^2)'$ be the parameter and by $\Theta$ be the parameter space.
 Assume that the observations $\{y_0, y_1,..., y_n\}$ are
from model (\ref{eq.model}) with the true value $\theta_0$.
Clearly, under Assumption \ref{assumption} below, it follows that
$\E(y_t|y_{t-1})=p\phi |y_{t-1}|$ and $\Var(y_t|y_{t-1})=p(1-p)\phi^2y_{t-1}^2+\sigma^2$.
Then the (conditional) log-quasi-likelihood function (omitting a constant) is
\begin{flalign*} 
L_n(\theta)=\sum_{t=1}^n\ell_t(\theta):=
\sum_{t=1}^n\left\{
\log\left[p(1-p)\phi^2y_{t-1}^2+\sigma^2\right]
+\frac{(y_t-p\phi |y_{t-1}|)^2}
{p(1-p)\phi^2y_{t-1}^2+\sigma^2}
\right\}.
\end{flalign*}
The QMLE of $\theta_0$ is defined as
\begin{flalign*}
\widehat{\theta}_n=\arg\min_{\theta\in\Theta}L_n(\theta).
\end{flalign*}

To study the asymptotics of $\widehat{\theta}_n$, the following two assumptions are needed.
\begin{assumption}\label{assumption}
$\{\varepsilon_t\}$ is i.i.d. and independent of i.i.d. binary variables $\{s_t\}$ with $p_0<1$.
Further, $\varepsilon_1$ has a positive density on $\mathbb{R}$ with zero mean and finite variance.
\end{assumption}
\begin{assumption}\label{para:space}
The parameter space $\Theta$ is a compact subset of $\{\theta=(\phi, p, \sigma^2)': \phi\neq0,\,
0< p<1, \, 0<\sigma^2<\infty\}$.
\end{assumption}

The following theorems states the strong consistency and the asymptotic normality of $\widehat{\theta}_n$.
\begin{theorem}\label{thm.consistent}
If Assumptions $\ref{assumption}$-$\ref{para:space}$ hold,
then $\widehat{\theta}_n\rightarrow \theta_0$ a.s.
as $n\rightarrow\infty$.
\end{theorem}

\begin{theorem}\label{thm.asym}
If Assumptions $\ref{assumption}$-$\ref{para:space}$ hold,
$\mathbb{E}(\varepsilon_t^4)<\infty$, and $\theta_0$ is
an interior point of $\Theta$, then
\begin{flalign*}
\sqrt{n}\big(\widehat{\theta}_n-\theta_0\big)\stackrel{d}\longrightarrow
\mathcal{N}\left(0, \,\mathcal{J}^{-1}\mathcal{I}\mathcal{J}^{-1}\right),\quad \mbox{as $n\to\infty$},
\end{flalign*}
where `$\stackrel{d}\longrightarrow$' stands for convergence in distribution,
\begin{flalign*}
\mathcal{J}
=\mathbb{E}\left\{\frac{1}{[p_0(1-p_0)\phi^2_0y_{t}^2+\sigma^2_0]^2}\mathbf{A}_t\right\}
+\mathbb{E}\left\{\frac{2y_{t}^2}{p_0(1-p_0)\phi^2_0y_{t}^2+\sigma^2_0}\right\}
\mathbf{D}, &&
\end{flalign*}
\begin{flalign*}
\mathcal{I}
&= \mathbb{E}\left\{\frac{p_0(1-p_0)(1-2p_0)^2\phi^4_0y_{t}^4
+4\sigma_0^2p_0(1-p_0)\phi_0^2y_{t}^2+(\kappa_4-\sigma^4_0)}{[p_0(1-p_0)\phi^2_0y_{t}^2+\sigma^2_0]^4}
\mathbf{A}_t\right\}\\
&\quad+\mathbb{E}\left\{\frac{4y^2_{t}}{p_0(1-p_0)\phi^2_0y_{t}^2+\sigma^2_0}\right\}\mathbf{D}
+\mathbb{E}\left\{\frac{2p_0(1-p_0)(1-2p_0)\phi_0^3y_{t}^4+2|y_{t}|\kappa_3}
{[p_0(1-p_0)\phi^2_0y_{t}^2+\sigma^2_0]^3}\mathbf{B}_t\right\}&&
\end{flalign*}
with $\kappa_3=\mathbb{E}(\varepsilon_t^3)$, $\kappa_4=\mathbb{E}(\varepsilon_t^4)$, and
\begin{flalign*}
\mathbf{A}_t&=
\left(
    \begin{array}{ccc}
      4p^2_0(1-p_0)^2\phi^2_0 y_{t}^4 & 2p_0(1-p_0)(1-2p_0)\phi^3_0 y_{t}^4 & 2p_0(1-p_0)\phi_0 y_{t}^2 \\
       & (1-2p_0)^2\phi^4_0 y_{t}^4 & (1-2p_0)\phi^2_0 y_{t}^2 \\
       &  & 1
    \end{array}
  \right),\\
\mathbf{B}_t&=\left(
   \begin{array}{ccc}
     4p_0^2(1-p_0)\phi_0 y_{t}^2 & p_0(3-4p_0)\phi^2_0 y_{t}^2 & p_0 \\
      & 2(1-2p_0)\phi^3_0 y_{t}^2 & \phi_0 \\
     &  & 0
   \end{array}
 \right),\quad \mbox{and}\quad
\mathbf{D}=\left(
  \begin{array}{ccc}
    p^2_0 & p_0\phi_0 & 0 \\
     & \phi^2_0 & 0 \\
     &  & 0
  \end{array}
\right).&&
\end{flalign*}
Here, the elements in the lower triangles can be completed by symmetry.
\end{theorem}

\textbf{Remark 1}. From the expressions of the matrices $\mathcal{I}$ and $\mathcal{J}$
in Theorem \ref{thm.asym}, we can see that
each element of random matrices within the expectation
is bounded and thus it is unnecessary to require moment conditions on $y_t$ for
the asymptotics of $\widehat{\theta}_n$.

\textbf{Remark 2}. In practice, to make statistical inference on $\theta_0$, we must estimate the matrices
$\mathcal{I}$ and $\mathcal{J}$. From the proof of Theorem \ref{thm.asym}, they can be consistently estimated by
\begin{flalign*}
\widehat{\mathcal{I}}_n=\frac{1}{n}\sum_{t=1}^n\frac{\partial\ell_t(\widehat{\theta}_n)}{\partial\theta}
\frac{\partial\ell_t(\widehat{\theta}_n)}{\partial\theta'}\quad\mbox{and}\quad
\widehat{\mathcal{J}}_n=\frac{1}{n}\sum_{t=1}^n\frac{\partial^2\ell_t(\widehat{\theta}_n)}
{\partial\theta\partial\theta'},
\end{flalign*}
respectively. Note that  the plug-in method is here invalid
since both $\kappa_3$ and $\kappa_4$ in $\mathcal{I}$
cannot be estimated from the residuals.
Additionally, due to the constraint $p_0\in (0, 1)$,  the Delta method
may be needed to construct confidence intervals of $p_0$. If necessary, for example,
we can consider the transformation
$g(p)=\log[(1-p)/p]$ for $p\in (0, 1)$. Note that
\begin{flalign*}
\sqrt{n}\big(g(\widehat{p}_n)-g(p_0)\big)
\stackrel{d}{\longrightarrow}\frac{\lambda_p}{p_0(1-p_0)}\mathcal{N}(0, 1),
\end{flalign*}
if $\sqrt{n}(\widehat{p}_n-p_0)\stackrel{d}{\longrightarrow}
\mathcal{N}(0, \lambda_p^2)$. Then, for any fixed $\alpha\in(0, 1)$,
a $100(1-\alpha)\%$ confidence interval of $p_0$ is
\begin{flalign*}
\left[\Big\{1+\exp\Big[g(\widehat{p}_n)
-\frac{\lambda_p\,z_{\alpha/2}}{\sqrt{n}\widehat{p}_n(1-\widehat{p}_n)}\Big]\Big\}^{-1},\,\:
\Big\{1+\exp\Big[g(\widehat{p}_n)
+\frac{\lambda_p\,z_{\alpha/2}}{\sqrt{n}\widehat{p}_n(1-\widehat{p}_n)}\Big]\Big\}^{-1}\right],
\end{flalign*}
where $z_{\alpha/2}$ is the lower $\alpha/2$-quantile of the standard normal.

\section{Model Diagnostic Checking}\label{section.checking}
Diagnostic checking is important for time series modeling. The
most commonly used tool is the portmanteau test, which depends on the autocorrelation of the
residuals or the squared residuals, see, e.g., \cite{McLeod1983}, \cite{Li1994}, \cite{LiWK2004},
 and \cite{Chen:Zhu2015}. However,
such the portmanteau test fails for the adequacy of model \eqref{eq.model} since
the residuals cannot be obtained. In fact, the residuals should be theoretically calculated by
$\widehat{\varepsilon}_t=y_t-s_t\widehat{\phi}_n|y_{t-1}|$ with the initial value $y_0$ for $i=1,...,n$.
Unfortunately, the latent variables $\{s_t: 1\leq t\leq n\}$ are unknown and prevent us from getting $\{\widehat{\varepsilon}_t\}$.

To check the adequacy of model \eqref{eq.model}, we introduce a new portmanteau test, which is constructed via a transformation of an uncorrelated sequence.
Note that the sequence $\{y_t-p_0\phi_0|y_{t-1}|: t\in\mathbb{Z}\}$ is still uncorrelated when $\mathbb{E}y_t^2<\infty$ after replacing
$s_t$ by its mean $p_0$ in \eqref{eq.model}. To reduce the dependence on the moments of $y_t$,
similar to \cite{Ling:2005, Ling:2007}, we adopt a self-weight method and then define a new sequence
$\{\eta_t\}$ by
\begin{flalign}\label{eta:expression}
\begin{split}
\eta_t:=&\eta_{t,a}=(y_t-p_0 \phi_0 |y_{t-1}|)I(|y_{t-1}|\leq a)\\
=&\phi_0(s_t-p_0)|y_{t-1}|I(|y_{t-1}|\leq a)+\varepsilon_tI(|y_{t-1}|\leq a),\quad
t\in\mathbb{Z},
\end{split}
\end{flalign}
where the constant $a$ is positive and is called a tuning parameter, and $I(\cdot)$ is an indicator function. Clearly, $\{\eta_{t}\}$ is
strictly stationary and ergodic  since it is a measurable function of strictly stationary and ergodic sequence
$(y_{t-1}, s_t, \varepsilon_t)'$.
Further, by the mutually independence among $s_t$, $\varepsilon_t$, and $y_{t-1}$,
a simple calculation yields that
\begin{flalign}\label{eta_sigma}
\E\eta_{t}=0,\,\,
\sigma_\eta^2:=\E\eta_{t}^2=p_0(1-p_0)\phi_0^2\E\left\{y_1^2I(|y_{1}|\leq a)\right\}
+\sigma_0^2\mathbb{P}(|y_{1}|\leq a), \,\, \E(\eta_{t}\eta_{t-k})=0,
\end{flalign}
for $k\geq 1$. That is, $\{\eta_t\}$ is always a white noise under Assumption \ref{assumption}.
Moreover, it is also a martingale difference sequence.

Let $\widehat\eta_{t}=(y_t-\widehat p_n \widehat \phi_n |y_{t-1}|)I(|y_{t-1}|\leq a)$, $1\leq t\leq n$. Intuitively,
its sample autocorrelation $\widehat \rho_{nk}$ should be close to zero if model specification is correct,
where
\begin{flalign*}
\widehat \rho_{nk}=\frac{\sum_{t=k+1}^n(\widehat \eta_{t}-\bar \eta)(\widehat \eta_{t-k}-\bar \eta)}{\sum_{t=1}^n(\widehat \eta_{t}-\bar \eta)^2}\quad\mbox{with}\quad
\bar \eta=n^{-1}\sum_{t=1}^n\widehat \eta_{t}.
\end{flalign*}
Denote $\widehat{\bm\rho}_{n}=(\widehat \rho_{n1},\dots,\widehat \rho_{nM})'$,
where $M\ge 1$ is a fixed positive integer. The following theorem gives
the limiting distribution of $\widehat{\bm\rho}_n$.

\begin{theorem}\label{thm.rho}
Suppose the conditions in Theorem $\ref{thm.asym}$ hold.
If model \eqref{eq.model} is correctly specified, then $\sqrt{n} \widehat{\bm\rho}_n\stackrel{d}\longrightarrow\
\mathcal{N}(0,\,\mathbf{UGU}')$, where $\mathbf{G}=\E [v_tv_t']$ with
\begin{flalign*}
v_t=\left(\frac{\eta_t\eta_{t-1}}{\sigma_\eta^2},\dots,
\frac{\eta_t\eta_{t-M}}{\sigma_\eta^2},
\left(-\mathcal{J}^{-1}
\frac{\partial \ell_t(\theta_0)}{\partial\theta}\right)'\,\right)',
\end{flalign*}
and $\mathbf{U}=[\mathbf{I}_M,\sigma_\eta^{-2}(u_{1},\dots, u_{M})'(p_0,\phi_0,0)]$
with $u_k=-\E\{\eta_{t-k}|y_{t-1}|I(|y_{t-1}|\leq a)\}$ and
$\sigma_\eta^2$ being defined in $(\ref{eta_sigma})$.
\end{theorem}

Based on Theorem \ref{thm.rho}, our portmanteau test statistic is defined as
\begin{flalign*}
Q_M=n\widehat{\bm\rho}_n^{\,'} \big(\widehat{\mathbf{U}}_n \widehat{\mathbf{G}}_n \widehat{\mathbf{U}}_n'\big)^{-1} \widehat{\bm\rho}_n,
\end{flalign*}
where $\widehat{\mathbf{U}}_n$ and $\widehat{\mathbf{G}}_n$ are consistently sample counterparts of
$\mathbf{U}$ and $\mathbf{G}$, respectively. Under conditions of Theorem \ref{thm.rho}, we have that
$Q_M \stackrel{d}\longrightarrow\ \chi^2_M$.

\textbf{Remark 3}. In application, we must choose the tuning parameter $a$
when our test statistic $Q_M$ is used.
According to the suggestion in the literature, see, for example, \cite{Ling:2005, Ling:2007},
we can let $a$ be the 90\% or 95\% quantile of data $\{|y_1|,...,|y_n|\}$. Many practical experience
shows that this self weight performs well,  although it may not be optimal and there exist some other choices. Further,
from Section \ref{stable}, we can see that $\E y_t^4<\infty$ if $p_0\phi_0^4<1$. Note that the condition
$p_0\phi_0^4<1$ is verifiable or testable. Thus, if $\widehat{p}_n\widehat{\phi}_n^{\,4}<1$, then we can let $a=+\infty$, that is, the self-weight is redundant.

\section{Bubble Tagging}\label{section.label}
 An important problem in economic or financial data analysis is to tag the bubbles, which can help the government or financial institutions to respond timely to resolve a potential financial crisis. Being able to successfully tag a bubble can also create lucrative trading opportunities. On the other hand, many economic studies can benefit from meaningful tagging of history bubbles to help understand the economic status and explain the reasoning behind certain economic behaviors in the history. The problem of bubble tagging, however, is often not easy and requires sophisticated statistical modeling and treatment. For this, \cite{Phillips2011Dating} considered decomposing the asset price process into a fundamental component determined by expected future dividends and an explosive bubble component, and proposed a recursive testing procedure. \cite{Phillips2015} modeled the null hypothesis as a random walk with asymptotically negligible drift and studied the limit theory of a dating algorithm for bubble detection; see also \citet{Phillips2015Testing}. We also refer to recent papers by \cite{Blasques2022A} and \cite{Kurozumi}, and references therein, for additional literature. Existing results on bubble detection, nevertheless, were mainly developed in a nonstationary framework for which the bubble mechanism is not incorporated into the underlying stationary process and treated as a separated period on the timeline.

A distinguishable feature of the current paper is to incorporate the bubble mechanism into the data generating process to provide a stationary statistical model that can capture and interpret bubbles. Unlike existing results where bubbles are assumed to persist for an adequate duration to achieve consistent detection, bubble tagging in the current stationary framework can be more challenging as bubbles, especially transient bubbles or bubbles that only last for a very short time, can be easily mixed with large white noise observations. We in the following provide two different methods for bubble tagging in the current stationary framework.

\subsection{A Residual-Based Method for Bubble Tagging}\label{subsec:bubbletaggingresidualbased}
Given model (\ref{eq.model}), we consider the difference
\begin{flalign}\label{eq.res}
r_t:=y_t-\phi_0|y_{t-1}|
=\left \{
\begin{array}{ll}
\varepsilon_t, & \mbox{if $s_t=1$},\\
\varepsilon_t-\phi_0|y_{t-1}|, &\mbox{if $s_t=0$}.
\end{array}
\right.
\end{flalign}
Since $\phi_0>1$ is generally assumed in real applications, it is expected that $r_t$ will be relatively smaller for time points where $s_t = 0$. Therefore, a natural approach is to tag $s_t = 0$ if the difference $r_t \leq c_r$ for some threshold $c_r$. We in the following provide some theoretical understanding of such a tagging method. To this end, we introduce an auxiliary process $\{z_t\}$, where $z_0 = \epsilon_0$ and
\begin{flalign}\label{eqn:zt}
z_t = \phi_0|z_{t-1}| + \varepsilon_t,\quad t = 1,2,\ldots.
\end{flalign}
Unlike the full model specified in (\ref{eq.model}) that is stationary for which a bubble can collapse, the process $\{z_t\}$ defined above is a pure bubble process that is explosive and nonstationary. In particular, for any $k \geq 1$, $z_k$ shares the same distribution as $y_{t+k}$ if a bubble forms at time $t+1$ and persists through time $t+k$, which we call a $k$-th cumulative bubble. This also relates to the excursive period with duration $k$ in financial applications; see for example our data analysis in Section \ref{section.example}. For consistent tagging of bubbles, it is generally required that $k \to \infty$, namely the bubble has to persist for a growing horizon of time; see for example \cite{Phillips2015} and references therein.

\begin{proposition}\label{prop:kbubblenbm}
For any time $t$, if the innovation distribution is symmetric, then the conditional probability that the collapse of a $k$-th cumulative bubble will be correctly tagged by the aforementioned method equals to $\mathbb{P}(z_k \geq -c_r)$, namely the marginal probability that the auxiliary explosive bubble process will exceed the same threshold in the other direction.
\end{proposition}

In practice, a threshold of $c_r < 0$ is typically chosen for $r_t$ defined in (\ref{eq.res}), and as a result $-c_r > 0$ will be a positive threshold for $z_k$. Given the explosive nature of the bubble process $\{z_k\}$, it is expected that $\mathbb{P}(z_k > -c_r) \to 1$ as $k \to \infty$ for any chosen threshold $-c_r > 0$, and as a result the probability that the collapse of a $k$-th cumulative bubble will be correctly tagged increases to one as $k \to \infty$. This resonates the result of \citet{Phillips2015} but in very different settings. To be more specific, \citet{Phillips2015} assumed that the bubble period is a deterministic segment with an increasing number of time points within the whole observation period, while the current setting treats the bubble as an integrated part of an underlying stationary process in (\ref{eq.model}).

In real applications the parameter $\phi_0$ in (\ref{eq.res}) is unknown and we propose to plug in the QMLE and tag $\widehat s_t = 0$ if the residual $\widehat r_t=y_t-\widehat \phi_n|y_{t-1}| < c_r$ for some threshold $c_r$. We in the following provide some empirical reference rules for choosing the threshold $c_r$.
\begin{itemize}
\item Rule 1 (hard threshold). Set $c_r$ as the $(1-\widehat p_n)$-th quantile of $\{\widehat r_t\}$. Such a choice is simple yet effective, and it can be seen from our simulation results in Section \ref{section.simulation} that it is also reasonably robust to different choices of innovation distributions.
\end{itemize}
If in particular the innovations are normal with $\varepsilon_t\sim \mathcal{N}(0, \sigma_0^2)$, then we can in addition consider the following choice of thresholds. For this, we use $(\widehat{\phi}_n,\widehat{p}_n,\widehat{\sigma}^2_n)$ to denote the QMLE of $(\phi_0,p_0,\sigma_0^2)$ as in Section \ref{section.qmle}.
\begin{itemize}
\item Rule 2 (conditional likelihood). From (\ref{eq.res}), conditioning on $y_{t-1}$, $r_t=\varepsilon_t+\mu$, where $\mu$ can be considered as a location parameter and takes only two values, i.e., $\mu=0$ or $-\phi_0|y_{t-1}|$, corresponding to $s_t=1$ or $s_t=0$, respectively. By comparing conditional likelihood of $r_t$ given $y_{t-1}$ for each time $t$, we can determine the value $\mu=0$ or $-\phi_0|y_{t-1}|$ and then determine whether we need to tag this time $t$.  Equivalently, for each time $t$, we can set $c_{r,t}=-\widehat \phi_n |y_{t-1}|/2$ and
    tag the time $t$ if $\widehat r_t < c_{r,t}$.
\item Rule 3 (time-varying quantile). Conditioning on $y_{t-1}$, we have $r_t\sim p_0 \mathcal N(0, \sigma_0^2)+(1-p_0)\mathcal N(-\phi_0|y_{t-1}|, \sigma_0^2)$.
Based on such a conditional distribution, a natural choice for $c_{r,t}$ is
\begin{flalign*}
c_{r,t}=\inf\left\{r\in\mathbb{R}:
\widehat p_n\Phi\left(r/\widehat{\sigma}_n\right)+
(1-\widehat{p}_n)\Phi\big((r+\widehat{\phi}_n|y_{t-1}|)/\widehat{\sigma}_n\big)>1-\widehat{p}_n
\right\}
\end{flalign*}
for each time $t$, where $\Phi(\cdot)$ is the cdf of $\mathcal N(0,1)$.
\item Rule 4 (Bayesian). Based on the Bayes' rule, conditioning on $y_{t-1}$, we can obtain the posterior probability mass function of $s_t$ given $r_t$, that is,
$\mathbb{P}(s_t=1|r_t)$ and $\mathbb{P}(s_t=0|r_t)$. Theoretically, for each time $t$, if
$\mathbb{P}(s_t=1|r_t)<\mathbb{P}(s_t=0|r_t)$, then, we can tag this time $t$. Otherwise, we do not do it. Equivalently, we can tag the time $t$ if
\begin{flalign*}
(1-p_0)f((r_t+\phi_0|y_{t-1}|)/\sigma_0)<p_0 f(r_t/\sigma_0),
\end{flalign*}
where $f(\cdot)$ is the density of $\mathcal{N}(0, 1)$. Using $\{\widehat{r}_t\}$,
we can tag the time $t$ if
\begin{flalign*}
(1-\widehat{p}_n) f\big(\big(\widehat{r}_t+\widehat{\phi}_n|y_{t-1}|\big)/\widehat{\sigma}_n\big)<
\widehat{p}_n f\left(\widehat{r}_t/\widehat{\sigma}_n\right).
\end{flalign*}
\end{itemize}

\subsection{A Null-Based Method for Bubble Tagging}\label{subsec:bubbletaggingnullbased}
The method described in Section \ref{subsec:bubbletaggingresidualbased} relies on residuals from the one-step ahead recursion specified by model (\ref{eq.model}) to tag the collapse of bubbles. In essence, it treats the explosive bubble alternative as the default and aims at detecting the null of no bubble as an anomaly. We shall here consider its complement which sets the null of no bubble as the baseline and detects the formation of a bubble as an anomaly. To be more specific, when $s_t = 0$ and there is no bubble at time $t$, we have $y_t = \varepsilon_t$ which forms a stationary white noise sequence. When the bubble starts to form at time $t$, however, an explosive drift $\phi_0|y_{t-1}|$ will be cumulatively added to the otherwise white noise sequence during the whole bubble period making the observed $y_t$ to cumulatively deviate away from the baseline. Therefore, it becomes natural to tag time $t$ as a bubble if $y_t > c$ for some threshold $c$. In contrast to the approach in Section \ref{subsec:bubbletaggingresidualbased} which relies exclusively on model (\ref{eq.model}) to compute the residuals $\{\widehat r_t\}$, this null-based method directly works on the original observations $\{y_t\}$ and can be less model dependent. In addition, since $y_t$ is distributed as a white noise sequence under the null of no bubble, the threshold $c$ can be taken as a uniform constant, which can be a convenient feature that facilitates the decision rule visualization. It can also be more advantageous in situations when bubbles are not prevailing in the observation period. Let $\{z_t\}$ be the auxiliary process defined in (\ref{eqn:zt}), and we in the following provide some theoretical understanding of such a null-based bubble tagging method under the fixed horizon domain.

\begin{proposition}\label{prop:kbubble}
For any time $t$, the conditional probability that a $k$-th cumulative bubble will be correctly tagged by the null-based method equals to $\mathbb{P}(z_k > c)$, namely the marginal probability that the auxiliary explosive bubble process will exceed the same threshold.
\end{proposition}

For bubbles that persist for a growing horizon of time, by the explosive nature of the auxiliary bubble process it is expected that $\mathbb{P}(z_k > c) \to 1$ as $k \to \infty$ for any given threshold $c$, and as a result the aforementioned null-based bubble tagging method can identify such a persistent bubble with probability tending to one. \citet{Phillips2015} treated the bubble period as a fixed but unknown deterministic section of the whole observation time, and provided the consistency when the length of the bubble section grows proportionally with the sample size. In contrast, the current paper treats the bubble as an intrinsic feature of a stationary data generating mechanism, which serves as an important step to provide a statistical model to understand the mechanism of an economic phenomenon. We also remark that, unlike the QMLE discussed in Section \ref{section.qmle}, the aforementioned null-based bubble tagging method and Proposition \ref{prop:kbubble} will continue to hold for situations when the hidden state process $\{s_t\}$ exhibits dependence and forms a stationary or nonstationary process by itself. For example, it can be a stationary Markov chain or a nonstationary Markov chain with time-varying transition matrices. In addition, the proof of Proposition \ref{prop:kbubble} can be readily generalized to handle bubble mechanisms other than the one-step autoregressive recursion specified in (\ref{eq.model}).

\section{Simulation Studies}\label{section.simulation}
To assess the performance of the QMLE of $\theta_0$ and Rules 1-4 in finite samples, we use the sample size
$n=200$, 400, and 800, each with 1000 replications for model (\ref{eq.model}).
The error $\varepsilon_t$ follows
\begin{itemize}
\item $\mathcal{N}(0, 1)$;
\item the Laplace distribution with density
\begin{flalign*}
h(x)=\frac{1}{\sqrt{2}}\exp\big(-\sqrt{2}|x|\big),\quad x\in\mathbb{R};
\end{flalign*}
\item the standardized Student's $t_5$ $(\mathrm{st}_5)$ with density
\begin{flalign*}
h(x)=\frac{8}{3\pi\sqrt{3}}\Big(1+\frac{x^2}{3}\Big)^{-3},\quad x\in\mathbb{R}.
\end{flalign*}
\end{itemize}
Three different true values
of $\theta_0=(\phi_0, p_0, \sigma^2_0)'$ are used, respectively, i.e.,
\begin{itemize}
\item Case I: $\theta_0=(1, 0.9, 1)'$;
\item Case II: $\theta_0=(\sqrt{10/9}, 0.9, 1)'$;
\item Case III: $\theta_0=(1.2, 0.9, 1)'$.
\end{itemize}
For Case I, $y_t$ is weakly stationary since $p_0\phi_0^2<1$, while $y_t$ is an infinite-variance process
in Case III since $p_0\phi_0^2>1$. For Case II, $\theta_0$ is on the boundary, i.e., $p_0\phi_0^2=1$,
which is never considered in the literature.

Table \ref{table.simu} reports the bias, empirical standard deviation (ESD), and asymptotic standard deviation (ASD) of the QMLE $\widehat{\theta}_n$ for Cases I-III.
Here, the ASD of $\theta_0$ is simulated by extra time series of length 10,000, and 2,000 replications are used to reduce the estimated bias.
From the table, we can see that the QMLE performs well irrespective of infinite variance or heavy-tailedness issues. The biases are small and all the ESDs are close to the corresponding ASDs.
\begin{table}[!htbp]\tiny
  \centering
  \caption{Numerical simulation results.\label{table.simu}}
    \begin{tabular}{clrrrrrrrrrrr}
    \toprule
    $n$ &    & \multicolumn{11}{c}{$(\phi_0,p_0,\sigma_0^2)$} \\
\cmidrule{1-1}\cmidrule{3-13}       &    & \multicolumn{3}{c}{$(1,0.9,1)$} &    & \multicolumn{3}{c}{$(\sqrt{10/9},0.9,1)$} &    & \multicolumn{3}{c}{$(1.2,0.9,1)$} \\
\cmidrule{3-5}\cmidrule{7-9}\cmidrule{11-13}       &    & \multicolumn{1}{c}{$\widehat\phi_n$} & \multicolumn{1}{c}{$\widehat p_n$} & \multicolumn{1}{c}{$\widehat\sigma^2_n$} &    & \multicolumn{1}{c}{$\widehat\phi_n$} & \multicolumn{1}{c}{$\widehat p_n$} & \multicolumn{1}{c}{$\widehat\sigma^2_n$} &    & \multicolumn{1}{c}{$\widehat\phi_n$} & \multicolumn{1}{c}{$\widehat p_n$} & \multicolumn{1}{c}{$\widehat\sigma^2_n$} \\
    \midrule
       & \multicolumn{12}{l}{$\varepsilon_t\sim\mathcal N(0,1)$} \\
    \midrule
    \multirow{3}[2]{*}{200} & Bias & -0.0012  & -0.0042  & 0.0098  &    & 0.0026  & -0.0077  & -0.0059  &    & 0.0018  & -0.0038  & -0.0060  \\
       & ESD & 0.0431  & 0.0467  & 0.1673  &    & 0.0307  & 0.0376  & 0.2047  &    & 0.0185  & 0.0284  & 0.1993  \\
       & ASD & 0.0367  & 0.0438  & 0.1559  &    & 0.0274  & 0.0373  & 0.1652  &    & 0.0176  & 0.0297  & 0.1913  \\
    \midrule
    \multirow{3}[2]{*}{400} & Bias & -0.0002  & -0.0026  & 0.0044  &    & 0.0009  & -0.0030  & 0.0014  &    & 0.0007  & -0.0020  & 0.0056  \\
       & ESD & 0.0265  & 0.0321  & 0.1112  &    & 0.0193  & 0.0266  & 0.1160  &    & 0.0126  & 0.0220  & 0.1378  \\
       & ASD & 0.0259  & 0.0310  & 0.1102  &    & 0.0194  & 0.0263  & 0.1168  &    & 0.0124  & 0.0210  & 0.1353  \\
    \midrule
    \multirow{3}[2]{*}{800} & Bias & 0.0002  & -0.0010  & 0.0009  &    & 0.0009  & -0.0021  & -0.0006  &    & 0.0005  & -0.0004  & -0.0009  \\
       & ESD & 0.0192  & 0.0217  & 0.0774  &    & 0.0145  & 0.0193  & 0.0864  &    & 0.0087  & 0.0145  & 0.0979  \\
       & ASD & 0.0183  & 0.0219  & 0.0779  &    & 0.0137  & 0.0186  & 0.0826  &    & 0.0088  & 0.0148  & 0.0957  \\
    \midrule
       & \multicolumn{12}{l}{$\varepsilon_t\sim$ Laplace} \\
    \midrule
    \multirow{3}[2]{*}{200} & Bias & 0.0020  & -0.0068  & 0.0143  &    & 0.0036  & -0.0093  & -0.0067  &    & 0.0035  & -0.0045  & -0.0098  \\
       & ESD & 0.0465  & 0.0492  & 0.2341  &    & 0.0323  & 0.0400  & 0.2350  &    & 0.0205  & 0.0313  & 0.2738  \\
       & ASD & 0.0399  & 0.0463  & 0.2199  &    & 0.0292  & 0.0384  & 0.2299  &    & 0.0183  & 0.0301  & 0.2626  \\
    \midrule
    \multirow{3}[2]{*}{400} & Bias & 0.0015  & -0.0056  & -0.0070  &    & 0.0018  & -0.0037  & 0.0034  &    & 0.0006  & -0.0015  & 0.0004  \\
       & ESD & 0.0289  & 0.0326  & 0.1545  &    & 0.0227  & 0.0286  & 0.1580  &    & 0.0135  & 0.0201  & 0.1950  \\
       & ASD & 0.0282  & 0.0328  & 0.1555  &    & 0.0206  & 0.0272  & 0.1626  &    & 0.0129  & 0.0213  & 0.1857  \\
    \midrule
    \multirow{3}[2]{*}{800} & Bias & -0.0001  & -0.0007  & 0.0034  &    & 0.0015  & -0.0024  & 0.0000  &    & 0.0007  & -0.0014  & -0.0040  \\
       & ESD & 0.0214  & 0.0235  & 0.1116  &    & 0.0150  & 0.0201  & 0.1179  &    & 0.0093  & 0.0151  & 0.1314  \\
       & ASD & 0.0199  & 0.0232  & 0.1099  &    & 0.0146  & 0.0192  & 0.1150  &    & 0.0092  & 0.0151  & 0.1313  \\
    \midrule
       & \multicolumn{12}{l}{$\varepsilon_t\sim \text{st}_5$} \\
    \midrule
    \multirow{3}[2]{*}{200} & Bias & 0.0002  & -0.0062  & 0.0046  &    & 0.0031  & -0.0085  & -0.0101  &    & 0.0029  & -0.0042  & -0.0074  \\
       & ESD & 0.0541  & 0.0527  & 0.2512  &    & 0.0345  & 0.0416  & 0.2461  &    & 0.0209  & 0.0315  & 0.2948  \\
       & ASD & 0.0433  & 0.0487  & 0.2712  &    & 0.0309  & 0.0395  & 0.2824  &    & 0.0191  & 0.0304  & 0.3224  \\
    \midrule
    \multirow{3}[2]{*}{400} & Bias & 0.0021  & -0.0050  & 0.0016  &    & 0.0018  & -0.0043  & 0.0063  &    & 0.0019  & -0.0028  & -0.0071  \\
       & ESD & 0.0312  & 0.0354  & 0.1833  &    & 0.0234  & 0.0291  & 0.1880  &    & 0.0149  & 0.0215  & 0.1910  \\
       & ASD & 0.0306  & 0.0345  & 0.1918  &    & 0.0219  & 0.0279  & 0.1997  &    & 0.0135  & 0.0215  & 0.2280  \\
    \midrule
    \multirow{3}[2]{*}{800} & Bias & -0.0001  & -0.0014  & 0.0083  &    & 0.0004  & -0.0019  & 0.0019  &    & 0.0009  & -0.0013  & -0.0033  \\
       & ESD & 0.0224  & 0.0239  & 0.1353  &    & 0.0150  & 0.0195  & 0.1312  &    & 0.0091  & 0.0151  & 0.1523  \\
       & ASD & 0.0216  & 0.0244  & 0.1356  &    & 0.0155  & 0.0197  & 0.1412  &    & 0.0096  & 0.0152  & 0.1612  \\
    \bottomrule
    \end{tabular}%
\end{table}%

\begin{figure}[t]
  \centering
  \includegraphics[width=6in]{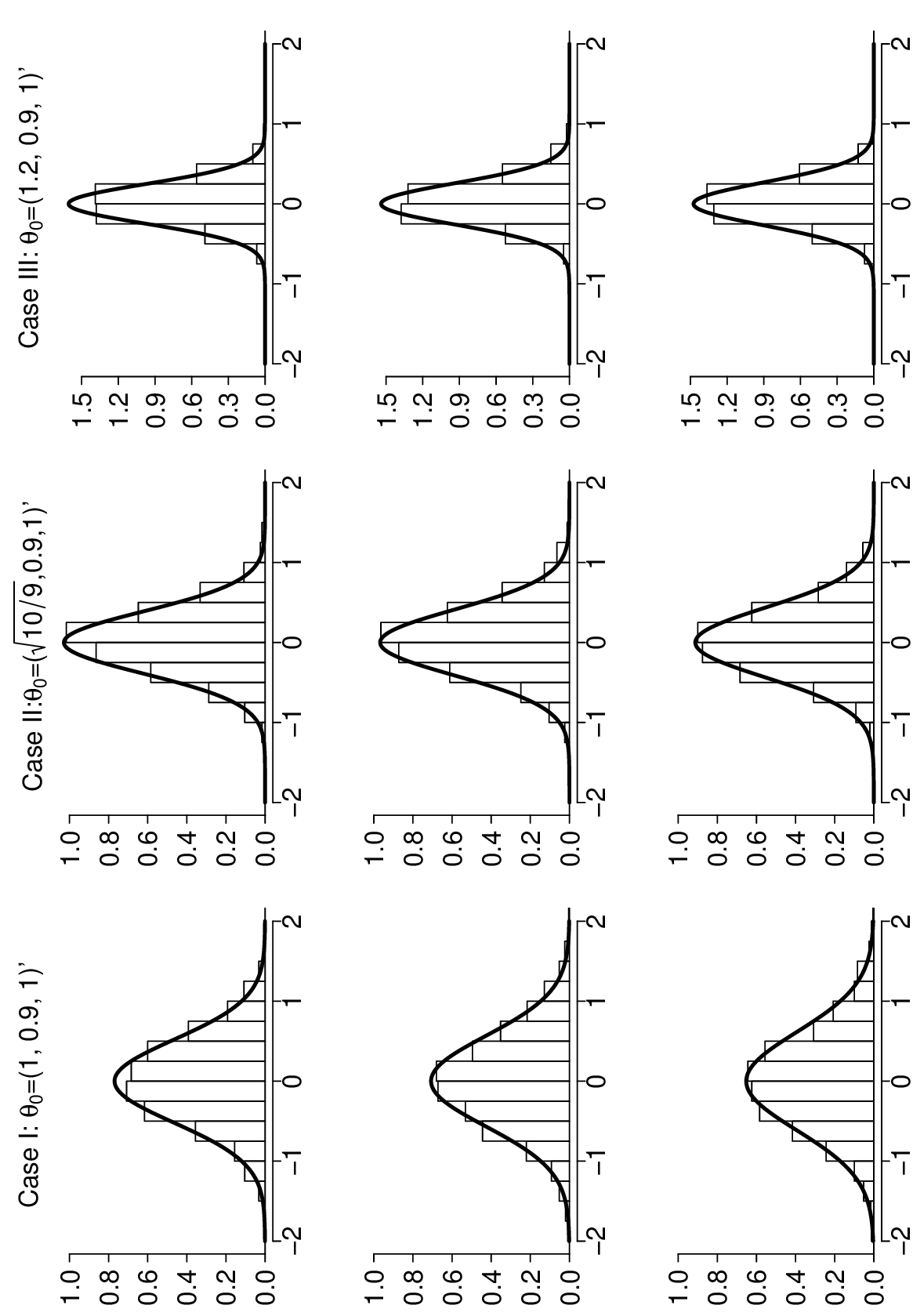}\\
  \caption{The histogram of $\sqrt{n}\big(\widehat{\phi}_n-\phi_0\big)$ with the sample size $n=400$.
The left column panel corresponds to Case I, i.e., $y_t$ is weakly stationary; the middle to Case II, and the right to Case III, i.e., $y_t$ has an infinite variance, respectively.
The upper row panel is when $\varepsilon_t\sim\mathcal{N}(0, 1)$, the middle when
$\varepsilon_t\sim$ the Laplace distribution, and the lower when $\varepsilon_t\sim \mathrm{st}_5$,
respectively.}\label{fig.AsymNormality}
\end{figure}

To see the overall approximation of the QMLE $\widehat{\phi}_n$, Fig.~\ref{fig.AsymNormality} displays
the histogram of $\sqrt{n}\big(\widehat{\phi}_n-\phi_0\big)$
when the sample size $n=400$. From the figure, we can see that $\sqrt{n}\big(\widehat{\phi}_n-\phi_0\big)$
is always asymptotically normal irrespective of infinite variance or heavy-tailedness of $y_t$.

We shall here examine the finite-sample performance of the two bubble tagging methods described in Section \ref{section.label}. For the residual-based tagging method in Section \ref{subsec:bubbletaggingresidualbased} with reference rules 1--4 we denote them by RBT$_1$--RBT$_4$ respectively in our numerical study, and we abbreviate the null-based tagging method in Section \ref{subsec:bubbletaggingnullbased} as NBT hereafter. For each generated process, let $\{\widehat s_t: 1\leq t\leq n\}$ be the estimated bubble tags and $\#$ denote the set cardinality. We consider the following evaluation metrics:
\begin{itemize}
\item P: the overall proportion of correct tagging $\#\{t:\widehat s_t=s_t\}/n$;
\item P0: the proportion of correctly tagged null states $\#\{t:\widehat s_t=0 , s_t=0\}/\#\{t:s_t=0\}$;
\item P1: the proportion of correctly tagged bubbles $\#\{t:\widehat s_t=1 , s_t=1\}/\#\{t:s_t=1\}$.
\end{itemize}
The results are summarized in Tables \ref{table.P.0.9} and \ref{table.P.0.5} based on 1000 replications for each configuration. To provide a fair comparison, we set the thresholds of different tagging methods so their estimated bubble ratios  $\#\{t:\widehat s_t=1\}/n$ are controlled at the same level. From Tables \ref{table.P.0.9} and \ref{table.P.0.5}, we can observe the followings.
\begin{enumerate}[(i)]
\item For both the RBT and NBT methods, the results are reasonably close across different error types. This indicates that the bubble tagging methods considered in Sections \ref{subsec:bubbletaggingresidualbased} and \ref{subsec:bubbletaggingnullbased} possess a certain degree of robustness with respect to the error distribution.
\item For each of the method considered, the performance in general improves when the nonlinear autoregressive coefficient $\phi_0$ increases. This is mainly because a larger value of the parameter $\phi_0$ in general leads to a stronger degree of explosiveness during the bubble period, making it relatively easier to distinguish between bubbles and null-states.
\item When $p_0 = 0.9$ as in Table \ref{table.P.0.9}, the performance of the RBT method can vary depending on which reference rule is used to obtain the threshold. The NBT method, on the other hand, seems to deliver a performance that is between the best and worst performed RBT methods. Note that the RBT method is deigned using the residuals that are more related to the bubble alternative, it meets with our intuition that the RBT method in general outperforms the NBT for most of the threshold choices when the bubble state probability $p_0 = 0.9$ is relatively high as in Table \ref{table.P.0.9}.
\item When the true underlying bubble state probability $p_0$ decreases to 0.5 as in Table \ref{table.P.0.5}, the bubble state no longer dominates and as a result the difference between the RBT and NBT methods becomes less noticeable and all the methods considered delivered quite similar performance.
\end{enumerate}

\begin{table}[htbp]
  \centering
  \caption{The values (in percentage) of P, P0, and P1 for RBT$_1$-RBT$_4$ and NBT when $n = 200$.}
    \begin{tabular}{crrrrrrrrrrr}
    \toprule
          & \multicolumn{3}{c}{$\phi_0=1, p_0=0.9$} &       & \multicolumn{3}{c}{$\phi_0=\sqrt{10/9}, p_0=0.9$} &       & \multicolumn{3}{c}{$\phi_0=1.2, p_0=0.9$} \\
\cmidrule{2-4}\cmidrule{6-8}\cmidrule{10-12}    Method & \multicolumn{1}{l}{P} & \multicolumn{1}{l}{P0} & \multicolumn{1}{l}{P1} &       & \multicolumn{1}{l}{P} & \multicolumn{1}{l}{P0} & \multicolumn{1}{l}{P1} &       & \multicolumn{1}{l}{P} & \multicolumn{1}{l}{P0} & \multicolumn{1}{l}{P1} \\
    \midrule
          & \multicolumn{11}{c}{$\varepsilon_t\sim\mathcal N(0,1)$} \\
    RBT$_1$ & 90.93  & 53.15  & 94.73  &       & 92.15  & 59.34  & 95.37  &       & 93.44  & 65.62  & 96.12  \\
    RBT$_2$ & 84.64  & 21.33  & 91.25  &       & 85.45  & 25.99  & 91.67  &       & 87.81  & 38.80  & 93.01  \\
    RBT$_3$ & 90.47  & 51.22  & 94.48  &       & 91.68  & 57.51  & 95.11  &       & 93.75  & 67.68  & 96.25  \\
    RBT$_4$ & 91.49  & 56.11  & 95.04  &       & 92.81  & 62.85  & 95.73  &       & 95.02  & 74.09  & 96.99  \\
    NBT   & 87.01  & 33.20  & 92.56  &       & 87.63  & 36.48  & 92.87  &       & 89.16  & 44.62  & 93.75  \\
          & \multicolumn{11}{c}{ $\varepsilon_t\sim$ Laplace } \\
    RBT$_1$ & 90.56  & 51.39  & 94.51  &       & 91.66  & 56.92  & 95.12  &       & 92.89  & 62.91  & 95.85  \\
    RBT$_2$ & 84.50  & 20.82  & 91.16  &       & 85.24  & 24.77  & 91.58  &       & 88.02  & 39.95  & 93.16  \\
    RBT$_3$ & 90.24  & 50.10  & 94.33  &       & 91.40  & 55.91  & 94.98  &       & 93.58  & 66.95  & 96.19  \\
    RBT$_4$ & 91.42  & 55.71  & 94.98  &       & 92.55  & 61.56  & 95.62  &       & 94.88  & 73.72  & 96.95  \\
    NBT   & 86.53  & 30.93  & 92.28  &       & 87.22  & 34.55  & 92.67  &       & 89.11  & 44.52  & 93.76  \\
          & \multicolumn{11}{c}{$\varepsilon_t\sim \text{st}_5$} \\
    RBT$_1$ & 90.63  & 51.80  & 94.56  &       & 91.87  & 58.01  & 95.20  &       & 93.10  & 63.99  & 95.94  \\
    RBT$_2$ & 84.49  & 20.64  & 91.16  &       & 85.33  & 25.62  & 91.59  &       & 87.97  & 39.76  & 93.10  \\
    RBT$_3$ & 90.36  & 50.61  & 94.41  &       & 91.61  & 57.11  & 95.06  &       & 93.77  & 67.86  & 96.26  \\
    RBT$_4$ & 91.40  & 55.71  & 94.98  &       & 92.81  & 62.86  & 95.72  &       & 94.97  & 74.05  & 96.97  \\
    NBT   & 86.69  & 31.66  & 92.38  &       & 87.47  & 35.97  & 92.77  &       & 89.22  & 44.92  & 93.79  \\
    \bottomrule
    \end{tabular}%
  \label{table.P.0.9}%
\end{table}%

\begin{table}[htbp]
  \centering
  \caption{The values (in percentage) of P, P0, and P1 for RBT$_1$-RBT$_4$ and NBT when $n = 200$.}
    \begin{tabular}{crrrrrrrrrrr}
    \toprule
          & \multicolumn{3}{c}{$\phi_0=1, p_0=0.5$} &       & \multicolumn{3}{c}{$\phi_0=\sqrt{10/9}, p_0=0.5$} &       & \multicolumn{3}{c}{$\phi_0=1.2, p_0=0.5$} \\
\cmidrule{2-4}\cmidrule{6-8}\cmidrule{10-12}    Method & \multicolumn{1}{l}{P} & \multicolumn{1}{l}{P0} & \multicolumn{1}{l}{P1} &       & \multicolumn{1}{l}{P} & \multicolumn{1}{l}{P0} & \multicolumn{1}{l}{P1} &       & \multicolumn{1}{l}{P} & \multicolumn{1}{l}{P0} & \multicolumn{1}{l}{P1} \\
    \midrule
          & \multicolumn{11}{c}{$\varepsilon_t\sim\mathcal N(0,1)$} \\
    RBT$_1$ & 66.70  & 66.79  & 66.69  &       & 67.77  & 67.87  & 67.75  &       & 69.93  & 70.05  & 69.92  \\
    RBT$_2$ & 68.13  & 68.24  & 68.11  &       & 69.45  & 69.56  & 69.43  &       & 72.16  & 72.32  & 72.15  \\
    RBT$_3$ & 68.09  & 68.20  & 68.08  &       & 69.37  & 69.48  & 69.35  &       & 71.93  & 72.09  & 71.92  \\
    RBT$_4$ & 68.32  & 68.43  & 68.31  &       & 69.55  & 69.66  & 69.53  &       & 72.42  & 72.58  & 72.41  \\
    NBT   & 66.85  & 66.95  & 66.83  &       & 67.97  & 68.07  & 67.95  &       & 70.00  & 70.16  & 69.99  \\
          & \multicolumn{11}{c}{ $\varepsilon_t\sim$ Laplace } \\
    RBT$_1$ & 67.84  & 68.01  & 67.98  &       & 68.75  & 68.85  & 68.73  &       & 70.60  & 70.72  & 70.63  \\
    RBT$_2$ & 70.07  & 70.26  & 70.21  &       & 71.28  & 71.41  & 71.26  &       & 73.55  & 73.70  & 73.58  \\
    RBT$_3$ & 69.95  & 70.15  & 70.09  &       & 71.25  & 71.38  & 71.24  &       & 73.53  & 73.67  & 73.55  \\
    RBT$_4$ & 70.18  & 70.38  & 70.32  &       & 71.39  & 71.52  & 71.37  &       & 73.81  & 73.97  & 73.83  \\
    NBT   & 68.13  & 68.32  & 68.28  &       & 69.00  & 69.12  & 68.99  &       & 70.75  & 70.88  & 70.77  \\
          & \multicolumn{11}{c}{$\varepsilon_t\sim \text{st}_5$} \\
    RBT$_1$ & 67.05  & 67.16  & 67.12  &       & 67.92  & 67.97  & 67.78  &       & 69.83  & 69.91  & 69.77  \\
    RBT$_2$ & 68.93  & 69.07  & 69.00  &       & 70.16  & 70.23  & 70.03  &       & 72.74  & 72.85  & 72.68  \\
    RBT$_3$ & 68.88  & 69.02  & 68.95  &       & 70.08  & 70.14  & 69.94  &       & 72.55  & 72.66  & 72.49  \\
    RBT$_4$ & 68.99  & 69.13  & 69.06  &       & 70.20  & 70.28  & 70.07  &       & 72.92  & 73.04  & 72.85  \\
    NBT   & 67.33  & 67.47  & 67.40  &       & 68.30  & 68.37  & 68.16  &       & 70.31  & 70.42  & 70.24  \\
    \bottomrule
    \end{tabular}%
  \label{table.P.0.5}%
\end{table}%

\clearpage
\section{An empirical example}\label{section.example}
\begin{figure}[!htbp]
  \centering
  \includegraphics[width=6in]{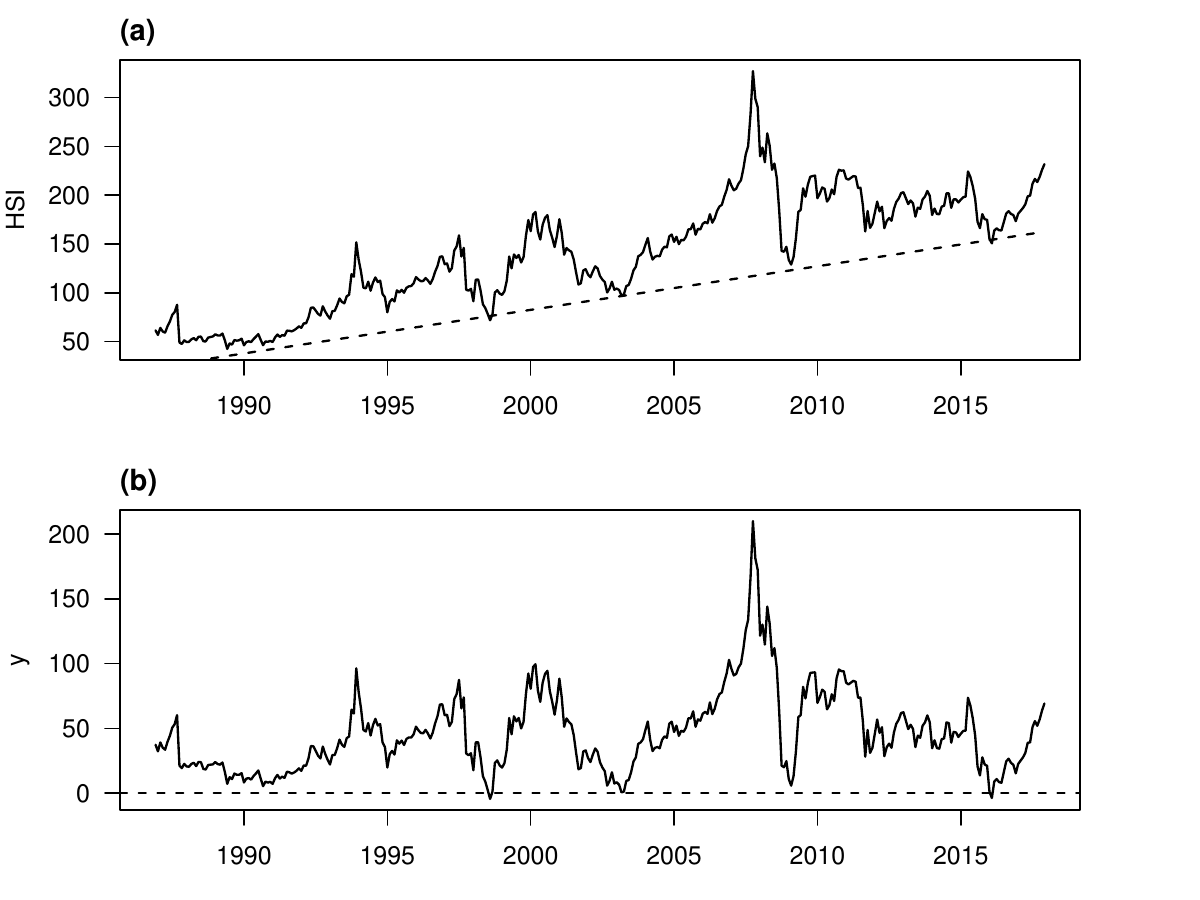}\\
  \caption{(a) Real HSI prices with the fitted linear trend (the dotted line); (b) $y_t$ series.}\label{fig.hsi}
\end{figure}

In this section, we analyze the monthly Hang Seng Index (HSI) from December 1986 to December 2017 with a total of 373 observations. To eliminate the effect of inflation on price, we transform nominal prices into real prices by the consumer price index, which can be obtained from the Federal Reserve Bank of St Louis.
Fig \ref{fig.hsi} (a) displays the real HSI prices, from which one can see an ascendant linear trend in the time series.
Thus, we first subtract such a linear trend from the series.
That is, we assume that the HSI real price $x_t$ is decomposed into
$$x_t=b_0+b_1 t +y_t,$$
where $b_0+b_1t$ denotes the linear trend and $y_t$ follows a SNAR model.
Note that $b_i,i=0,1$ can be seen as unknown parameters and can be estimated jointly.
Their estimates are $\widehat{b}_0=23.661$ and $\widehat{b}_1= 0.372$, respectively.
The linear time trend is plotted in Fig. \ref{fig.hsi} (a) by the dotted line and $\{y_t\}$
in Fig. \ref{fig.hsi} (b).
The estimates with standard deviations (SDs) of the SNAR model $\{y_t\}$ are reported in Table \ref{table.HSI}.
\begin{table}[!htbp]
  \centering
  \caption{The estimate with standard deviations (SD) of the fitted SNAR model.}
    \begin{tabular}{lrrr}
    \toprule
          & \multicolumn{1}{c}{$\phi_0$} & \multicolumn{1}{c}{$p_0$} & \multicolumn{1}{c}{$\sigma_0^2$} \\
    \midrule
    Estimate & 1.026  & 0.977  & 36.314  \\
    SD   & 0.011  & 0.005  & 8.649  \\
    \bottomrule
    \end{tabular}%
  \label{table.HSI}%
\end{table}%
All estimates are statistically significant at the $5\%$ level since their corresponding $p$-values are
extremely small which are thus not reported in the table.
The estimate of $\phi_0$ is large than one, and its $95\%$ confidence interval is $(1.005, 1.047)$, conforming to the locally explosive behavior of series $y_t$. For the fitting adequacy,
we calculate the $p$-values of the test statistic $Q_M$ with $M=6, 12, 18$, and $24$ when the tuning parameter
$a$ is the 90\% or 95\% quantile of $\{|y_t|,t=1\dots,n\}$, respectively. The results are summarized in Table
\ref{table.pvalues}, which implies that the fitting is adequate at the $5\%$ level.
\begin{table}[!htbp]
  \centering
  \caption{The $p$-values of $Q_M$. }
    \begin{tabular}{lcccc}
    \toprule
      $a\backslash M$    & 6     & 12    & 18    & 24 \\
    \midrule
    90\% & 0.7099  & 0.2427  & 0.3087  & 0.1549  \\
    95\% & 0.8588  & 0.7489  & 0.4898  & 0.1876  \\
    \bottomrule
    \end{tabular}%
  \label{table.pvalues}%
\end{table}%

\begin{figure}[t]
  \centering
  \includegraphics[width=15cm]{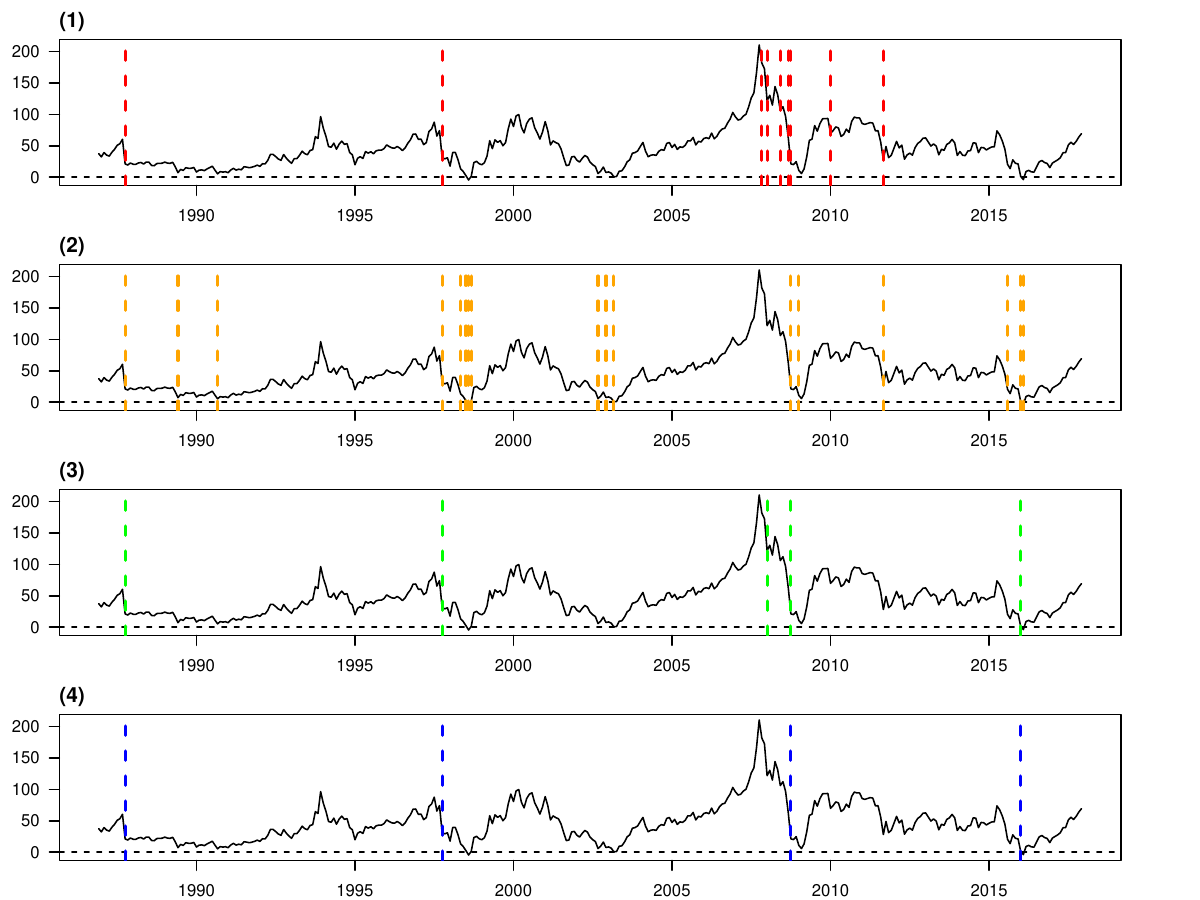}\\
  \caption{Selected dates of $\widehat{s}_t=0$ by Rules 1\textcolor{red}{--}4.}\label{fig.dates}
\end{figure}

We then apply the bubble tagging methods described in Section \ref{section.label} to label each time point as either being in a bubble state or being in the null. Since the estimated bubble probability $\widehat p_0 = 0.977$ from Table \ref{table.HSI} which is very high, in view of the simulation results in Section \ref{section.simulation}, we shall here consider using the residual-based method in Section \ref{subsec:bubbletaggingresidualbased} to tag the collapses of bubbles for the series $\{y_t\}$. In particular, Fig. \ref{fig.dates} displays the selected dates of $\widehat s_t=0$ under Rules 1--4. It can be seen from Fig. \ref{fig.dates} that the tagging times can vary based on which Rule is used, but several important dates are identified simultaneously by at least two rules. Table \ref{table.dates} summarizes such these dates,
\begin{table}[!htbp]
  \centering
  \caption{The selected important dates of $s_t=0$.}
    \begin{tabular}{lccccccc}
    \toprule
    Date & 1987-10  & 1997-10  & 2008-01 & 2008-10 & 2011-09 & 2016-01 \\
    Rule   & \{1,2,3,4\}  & \{1,2,3,4\}  & \{1, 3\} & \{1,2,3,4\} & \{1,2\} &\{2,3,4\}\\
    \bottomrule
    \end{tabular}%
  \label{table.dates}%
\end{table}%
which coinside with historical financial crises, i.e., the depression started from the Black Monday
in 1987, the Asian financial crises in 1997, the global financial turmoil caused
by the subprime crisis over 2007-2009, and the Hong Kong stock market plummeting in 2016.

\begin{figure}[!htbp]
  \centering
  \includegraphics[width=6in]{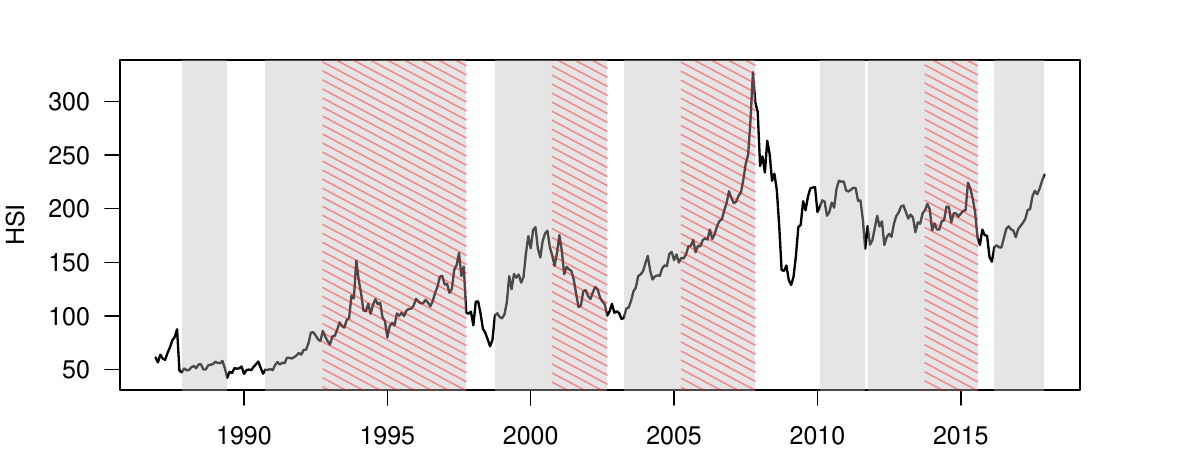}\\
  \caption{Excursive periods with duration exceeding 18 months shown by gray shadows; bubbles last over 24 months plotted by the shadow with red backslash.}\label{fig.period}
\end{figure}

Although the collapse of a bubble can be dated by $\widehat s_t=0$,
the emergence and exuberance of a bubble can not be asserted by $\widehat s_t=1$ immediately.
After all, a short-period deviation of the price is reasonable due to the market fluctuations.
Of course, a short-period deviation might be regarded as a small bubble in some sense, which bursts quickly by the market adjustment, thus we could pay little attention and ignore them afterwards.
What we really need to worry about is the bubble that can trigger tremendous harm, which emerges as the accumulation of long-lasting excursion.
Specifically, if $s_t=0, s_{t+1}=s_{t+2}=\dots=s_{t+k-1}=1, s_{t+k}=0$, then we call it an excursive period that starts from $t+1$ and ends at $t+k$, and define its duration as $k$.
Within an excursive period, the presence of a bubble can be supported if the duration exceeds some time span, for example, one or two years.
For our application, the time span is set to be 18 months.
Table \ref{table.period} summarizes the periods whose durations exceed 18 months,
as well as their start and end dates.
\begin{table}[!htbp]
  \centering
  \caption{Excursive period with duration exceeding 18 months.}
    \begin{tabular}{ccc}
    \toprule
    start & end   & duration \\
    \midrule
    1987-11 & 1989-06 & 20 \\
    1990-10 & 1997-10 & 85 \\
    1998-10 & 2002-09 & 48 \\
    2003-04 & 2007-11 & 56 \\
    2010-02 & 2011-09 & 20 \\
    2011-10 & 2015-08 & 47 \\
    2016-03 & 2017-12 & 22 \\
    \bottomrule
    \end{tabular}
  \label{table.period}
\end{table}
Fig. \ref{fig.period} plots those periods by gray shadows.
We can see explosive behaviors in most of the periods, indicating the presence and accumulation of bubbles. By Proposition \ref{prop:kbubble} in Section \ref{subsec:bubbletaggingresidualbased}, the RBT method is capable of detecting the collapse of an accumulated bubble when its duration $k \to \infty$; see also the same proposition for a probabilistic bound with a finite duration.
Another finding is that the magnitude of a bubble is larger as the period lasts longer possibly, for example, the one reaches a value of 210 in October 2007, corresponding to the period from April 2003 to November 2007 with the duration of 56 months.
Investors should be alert to such a long-time excursion along with the potential of disastrous bubbles.
In the periods where the bubble lasts over 24 months (plotted by the shadow with red backslash in Fig. \ref{fig.period}), one should be aware of the false boom in financial markets, and adjust asset allocation to hedge the risk of a potential bubble burst.

\section{Conclusions}\label{conclusion}
The paper has introduced a novel stochastic nonlinear autoregressive model to
describe the dynamics of economic or financial bubbles within a causal and stationary framework,
and discussed its strict stationarity and geometric ergodicity.
The paper has further studied the quasi-maximum likelihood estimation of the model
and established the asymptotics under minimal assumptions on innovation.
Due to the unobservability of the latent variable $s_t$ and
the resulting unavailability of the residuals, a new model diagnostic checking tool
has been proposed for the adequacy of the fitting. Finally, the paper considers two approaches, one from the residual perspective and the other from the null perspective, for bubble tagging.

Although our new model is useful, the model assumption on the independence between $\{\varepsilon_t\}$ and $\{s_t\}$
seems a little bit stronger from the perspective of empirical pragmatism. To obtain more reasonable interpretation or approximation of the bubble, such an
independence assumption can be relaxed. For  instance, we can assume that $s_t$ depends on the history of the observed process. Specifically, we can let $\mathbb{P}(s_t=1|\mathcal{F}_{t-1})=g(\bm\beta'\mathbf{y}_{t-1})$,
where $\mathcal{F}_{t-1}=\sigma(y_{t-j}: j\geq 1)$ be a sigma-field, $\mathbf{y}_{t-1}=(1, y_{t-1},...,y_{t-q})'$, and $g$ is a measurable function (e.g. a logistic function).
Furthermore, we can also restrict the form of $s_t$ in macroeconomic time series analysis
and let $s_t=I(\bm\beta'\mathbf{x}_t>c)$,
where $\mathbf{x}_t$ may contain many exogenous macroeconomic variables or indexes and $c$ is a threshold parameter. In addition, it is possible to consider the situation when the hidden state process $\{s_t\}$ exhibits temporal dependence and forms a Markov chain. In this case, the null-based bubble tagging method in Section \ref{subsec:bubbletaggingnullbased} can be more advantageous when bubbles occur in separated but persistent clusters. Another potential topic is to study multivariate stochastic nonlinear AR models.
We leave these topics for future research.



\appendix
\section{Technical Proofs}
\subsection{Proof of Theorem \ref{ssge}}
When $\phi_0=0$, then $\{y_t\}$ reduces to an i.i.d. sequence $\{\varepsilon_t\}$, and in this case all results hold clearly. Without loss of generality, we assume that $\phi_0\neq 0$ in what follows.
It suffices to verify the conditions in Theorem 19.1.3 in
\cite{MeynTweedie}. It is clear that $\{y_t\}$ defined by (\ref{eq.model}), with initial value $y_0$, is
an homogeneous Markov chain on $\mathbb{R}$ endowed with its Borel
$\sigma$-field $\mathcal{B}(\mathbb{R})$. Denote by $\lambda$ the Lebesgue measure on $(\mathbb{R}, \mathcal{B}(\mathbb{R}))$.
The transition probabilities of $\{y_t\}$ are given, for $y\in\mathbb{R}$,
$B\in \mathcal{B}(\mathbb{R})$, by
\begin{flalign*}
\mathbf{P}(y, B)=\mathbb{P}(y_t\in B|y_{t-1}=y)=p_0\mathbb{P}(\varepsilon_1+\phi_0|y|\in B)
+(1-p_0)\mathbb{P}(\varepsilon_1\in B).
\end{flalign*}

First, since $\mathbf{P}(\cdot, B)$ is continuous, for any $B\in \mathcal{B}(\mathbb{R})$, the chain $\{y_t\}$
has the Feller property.

Second, note that the density of $\varepsilon_1$ is positive over $\mathbb{R}$,
we have $\mathbf{P}(y, B)>0$ whenever $\lambda(B)>0$. Thus the chain $\{y_t\}$ is $\lambda$-irreducible.
Further, it can also be shown
that the $k$-step transition probabilities $\mathbf{P}^k(y, B)=\mathbb{P}(y_t\in B|y_{t-k}=y)=\int_{\mathbb{R}}\mathbf{P}^{k-1}(x, B)\mathbf{P}(y, dx)>0$
by an inductive approach for any integer $k\geq 1$,
whenever $\lambda(B)>0$, which establishes the aperiodicity of the chain $\{y_t\}$.

Third, let $V(x)=\log(1+|x|)$, $x\in\mathbb{R}$. Then, by a simple calculation, it follows that
\begin{flalign*}
\mathbb{E}\{V(y_t)|y_{t-1}=y\}
=(1-p_0)\mathbb{E}\log(1+|\varepsilon_1|)+p_0\mathbb{E}\log(1+|\phi_0|y|+\varepsilon_1|).
\end{flalign*}
Thus, we have that
\begin{flalign*}
\lim_{|y|\rightarrow\infty}\frac{\mathbb{E}\{V(y_t)|y_{t-1}=y\}}{V(y)}
&=\lim_{|y|\rightarrow\infty}\frac{(1-p_0)\mathbb{E}\log(1+|\varepsilon_1|)}{\log(1+|y|)}
+p_0\lim_{|y|\rightarrow\infty}\frac{\mathbb{E}\log(1+|\phi_0|y|+\varepsilon_1|)}{\log(1+|y|)}\\
&=0+p_0\lim_{|y|\rightarrow\infty}\left(\frac{\log|y|}{\log(1+|y|)}
+\frac{\mathbb{E}\log\Big(\frac{1}{|y|}+\big|\phi_0+\frac{\varepsilon_1}{|y|}\big|\Big)}{\log(1+|y|)}\right)\\
&=p_0.
\end{flalign*}
Since $p_0\in[0, 1)$, for fixed $\delta\in(0, 1-p_0)$, i.e., $p_0<1-\delta<1$, there exists a constant $M>0$ such that
\begin{flalign*}
\mathbb{E}\{V(y_t)|y_{t-1}=y\}\leq (1-\delta) V(y), \quad \mbox{when $|y|>M$}.
\end{flalign*}

To sum up the above arguments, by Theorem 19.1.3 in \cite{MeynTweedie},  there exists a geometrically ergodic solution to model (\ref{eq.model}). The solution is unique since $\mathbb{E}\log|s_t\phi_0|=-\infty$.
Thus, the results hold and then the proof is complete. \hfill $\square$

\subsection{Proof of Theorem \ref{thm.consistent}} Consider $\beta_n(\theta):=\{L_n(\theta)-L_n(\theta_0)\}/n$,
$\theta\in\Theta$. By the strong law of large numbers for stationary and ergodic sequences and
the inequality $\log x +x^{-1}-1\geq0$ for $x>0$, a conditional argument yields that
\begin{flalign*}
\beta_n(\theta)&=\frac{1}{n}\sum_{t=1}^n\left\{\log \frac{p(1-p)\phi^2y_{t-1}^2+\sigma^2}
{p_0(1-p_0)\phi_0^2y_{t-1}^2+\sigma_0^2}+\frac{(y_t-p\phi |y_{t-1}|)^2}{p(1-p)\phi^2y_{t-1}^2+\sigma^2}
-\frac{(y_t-p_0\phi_0 |y_{t-1}|)^2}{p_0(1-p_0)\phi_0^2y_{t-1}^2+\sigma_0^2}\right\}\\
&\stackrel{\mathrm{a.s.}}{\longrightarrow}
\mathbb{E}\left\{\log \frac{p(1-p)\phi^2y_{t-1}^2+\sigma^2}
{p_0(1-p_0)\phi_0^2y_{t-1}^2+\sigma_0^2}+\frac{(y_t-p\phi |y_{t-1}|)^2}{p(1-p)\phi^2y_{t-1}^2+\sigma^2}
-\frac{(y_t-p_0\phi_0 |y_{t-1}|)^2}{p_0(1-p_0)\phi_0^2y_{t-1}^2+\sigma_0^2}\right\}\\
&=\mathbb{E}\left\{\log \frac{p(1-p)\phi^2y_{t-1}^2+\sigma^2}
{p_0(1-p_0)\phi_0^2y_{t-1}^2+\sigma_0^2}+\frac{p_0(1-p_0)\phi_0^2y_{t-1}^2+\sigma_0^2}{p(1-p)\phi^2y_{t-1}^2+\sigma^2}
-1+\frac{(p\phi-p_0\phi_0)^2y_{t-1}^2}{p(1-p)\phi^2y_{t-1}^2+\sigma^2}\right\}\geq 0,
\end{flalign*}
where the equality holds if and only if
\begin{flalign*}
p(1-p)\phi^2y_{t-1}^2+\sigma^2=p_0(1-p_0)\phi_0^2y_{t-1}^2+\sigma_0^2\quad\mbox{and}\quad
(p\phi-p_0\phi_0)^2=0\quad \mbox{a.s.},
\end{flalign*}
equivalently, $\left\{p(1-p)\phi^2-p_0(1-p_0)\phi_0^2\right\}y_{t-1}^2=\sigma_0^2-\sigma^2$ a.s.
Then
\begin{flalign*}
p(1-p)\phi^2-p_0(1-p_0)\phi_0^2=0, \quad \sigma_0^2-\sigma^2=0.
\end{flalign*}
Combining with $(p\phi-p_0\phi_0)^2=0$, we have $\phi=\phi_0$, $p=p_0$ and $\sigma^2=\sigma_0^2$, i.e.,
$\theta=\theta_0$. The remainder of the proof can be completed by a standard compactness argument
and it is thus omitted. \hfill $\square$

\subsection{Proof of Theorem \ref{thm.asym}}
Let $q_t(\theta)=p(1-p)\phi^2y_{t-1}^2+\sigma^2$ and $q_t:=q_t(\theta_0)$. Then
the first- and second-order partial derivatives of $q_t(\theta)$ with respect to $\theta$ are respectively as follows
\begin{flalign}\label{eq.A.1}
\frac{\partial q_t(\theta)}{\partial\theta}
=\left(
\begin{array}{c}
2p(1-p)\phi y_{t-1}^2 \\
(1-2p)\phi^2 y_{t-1}^2 \\
1
\end{array}
\right),\quad
\frac{\partial^2 q_t(\theta)}{\partial\theta\partial\theta'}
=\left(
\begin{array}{ccc}
2p(1-p)y_{t-1}^2 & 2(1-2p)\phi y_{t-1}^2 & 0 \\
  & -2\phi^2 y_{t-1}^2 & 0 \\
 &  & 0
\end{array}
\right). &&
\end{flalign}
Using the notation $q_t(\theta)$, we have
\begin{flalign*}
\ell_t(\theta)=\log\left[q_t(\theta)\right]+\frac{(y_t-p\phi |y_{t-1}|)^2}{q_t(\theta)}.
\end{flalign*}
A simple calculation yields the first-order partial derivatives of
$\ell_t(\theta)$ with respect to $\theta$
\begin{flalign}\label{eq.A.2}
\frac{\partial \ell_t(\theta)}{\partial\theta}&=
\frac{1}{q_t(\theta)}\frac{\partial q_t(\theta)}{\partial\theta}
-\frac{2|y_{t-1}|(y_t-p\phi |y_{t-1}|)}{q_t(\theta)}\vartheta
-\frac{(y_t-p\phi |y_{t-1}|)^2}{\left[q_t(\theta)\right]^2}\frac{\partial q_t(\theta)}{\partial\theta},
\end{flalign}
where $\vartheta=(p, \phi, 0)'$,  and the second-order partial derivatives
\begin{flalign*}
\frac{\partial^2 \ell_t(\theta)}{\partial\theta\partial\theta'}&
=\left\{\frac{1}{q_t(\theta)}-\frac{(y_t-p\phi |y_{t-1}|)^2}{[q_t(\theta)]^2}\right\}\frac{\partial^2 q_t(\theta)}{\partial\theta\partial\theta'}
+\left\{\frac{2(y_t-p\phi |y_{t-1}|)^2}{[q_t(\theta)]^3}-\frac{1}{[q_t(\theta)]^2}\right\}
\frac{\partial q_t(\theta)}{\partial\theta}
\frac{\partial q_t(\theta)}{\partial\theta'}
\\
&\quad+\frac{2|y_{t-1}|(y_t-p\phi |y_{t-1}|)}{[q_t(\theta)]^2}
\left\{\frac{\partial q_t(\theta)}{\partial\theta}\vartheta'+\vartheta\frac{\partial q_t(\theta)}{\partial\theta'}\right\}+
\frac{2y_{t-1}^2}{q_t(\theta)}\vartheta\vartheta'\\
&\quad-\frac{2|y_{t-1}|(y_t-p\phi |y_{t-1}|)}{q_t(\theta)}
\begin{bmatrix}
0&1&0\\
1&0&0\\
0&0&0
\end{bmatrix}.
\end{flalign*}
By the Taylor expansion, by the definition of $\widehat{\theta}_n$, we have
\begin{flalign*}
0=\frac{1}{\sqrt{n}}\frac{\partial L_n(\widehat{\theta}_n)}{\partial\theta}=
\frac{1}{\sqrt{n}}\frac{\partial L_n(\theta_0)}{\partial\theta}
+\frac{1}{n}\frac{\partial^2L_n(\theta^*)}{\partial\theta\partial\theta'}
\sqrt{n}\big(\widehat{\theta}_n-\theta_0\big),
\end{flalign*}
where $\theta^*\in \Theta$ and satisfies $\|\theta^*-\theta_0\|\leq \|\widehat{\theta}_n-\theta_0\|$.
Note that the continuity of $\partial^2 \ell_t(\theta)/\partial\theta\partial\theta'$ in $\theta$
and the strong law of large numbers for stationary and ergodic sequences,
it is not hard to get
\begin{flalign*}
\frac{1}{n}\frac{\partial^2L_n(\theta^*)}{\partial\theta\partial\theta'}=
\frac{1}{n}\frac{\partial^2L_n(\theta_0)}{\partial\theta\partial\theta'}+o_p(1)=\mathcal{J}+o_p(1).
\end{flalign*}

Further, let $\mathcal{F}_t=\sigma(y_j: j\leq t)$ be the $\sigma$-algebra generated
by the random variables $\{y_j: j\leq t\}$. By the expressions in \eqref{eq.A.1} and \eqref{eq.A.2}, and the
following facts
\begin{flalign}\label{eq.A.3}
\begin{split}
\mathbb{E}[ y_t-p_0\phi_0 |y_{t-1}| |\mathcal{F}_{t-1}]
&=\mathbb{E}[(s_t-p_0)\phi_0|y_{t-1}|+\varepsilon_t|\mathcal{F}_{t-1}]=0,\\
\mathbb{E}[(y_t-p_0\phi_0 |y_{t-1}|)^2|\mathcal{F}_{t-1}]&=p_0(1-p_0)\phi^2_0y_{t-1}^2+\sigma^2_0=q_t,
\end{split}
\end{flalign}
we have that
\begin{flalign*}
\mathbb{E}\Big\{\frac{\partial\ell_t(\theta_0)}{\partial\theta}\Big|\mathcal{F}_{t-1}\Big\}=0,
\end{flalign*}
i.e., $\{\partial\ell_t(\theta_0)/\partial\theta\}$ is a martingale difference sequence with respect to
$\{\mathcal{F}_{t}\}$. Thus, by the martingale central limit theorem in \cite{Brown1971MartingaleCL}, it follows that
\begin{flalign*}
\frac{1}{\sqrt{n}}\frac{\partial L_n(\theta_0)}{\partial\theta}=
\frac{1}{\sqrt{n}}\sum_{t=1}^n\frac{\partial\ell_t(\theta_0)}{\partial\theta}\stackrel{d}{\longrightarrow}
\mathcal{N}(0, \mathcal{I}),
\end{flalign*}
where
\begin{flalign}\label{eq.A.I}
\mathcal{I}=\mathbb{E}\left\{\frac{\partial\ell_t(\theta_0)}{\partial\theta}
\frac{\partial\ell_t(\theta_0)}{\partial\theta'}\right\}.
\end{flalign}
Finally, we have
\begin{flalign}\label{eq.A.4}
\sqrt{n}\big(\widehat{\theta}_n-\theta_0\big)=-[\mathcal{J}+o_p(1)]^{-1}
\frac{1}{\sqrt{n}}\frac{\partial L_n(\theta_0)}{\partial\theta}\stackrel{d}{\longrightarrow}
\mathcal{N}(0, \mathcal{J}^{-1}\mathcal{I}\mathcal{J}^{-1}).
\end{flalign}

As for the explicit expressions of $\mathcal{I}$ and $\mathcal{J}$, by \eqref{eq.A.1}, \eqref{eq.A.3},
and the following facts
\begin{flalign*}
\mathbb{E}[(y_t-p_0\phi_0 |y_{t-1}|)^3|\mathcal{F}_{t-1}]&
=p_0(1-p_0)(1-2p_0)\phi_0^3|y_{t-1}|^3+\kappa_3,\\
\mathbb{E}[(y_t-p_0\phi_0 |y_{t-1}|)^4|\mathcal{F}_{t-1}]&
=p_0(1-p_0)(1-3p_0+3p_0^2)\phi_0^4y_{t-1}^4+6\sigma_0^2p_0(1-p_0)\phi_0^2y_{t-1}^2+\kappa_4,
\end{flalign*}
a tedious algebraic calculation can yield them and the detail is omitted.
The proof is complete. \hfill $\square$

\subsection{Proof of Theorem \ref{thm.rho}}
According to the definition of $\widehat\eta_{t}$, by Theorem \ref{thm.consistent}
and the strong law of large numbers for a stationary and ergodic sequence,
we first have the following facts,
as $n\rightarrow\infty$,
\begin{flalign*}
\bar{\eta}
&=\phi_0\frac{1}{n}\sum_{t=1}^n(s_t-p_0)|y_{t-1}|I(|y_{t-1}|\leq a)
+\frac{1}{n}\sum_{t=1}^n\varepsilon_t I(|y_{t-1}|\leq a)\\
&\qquad+(p_0\phi_0-\widehat{p}_n\widehat{\phi}_n)\frac{1}{n}\sum_{t=1}^n|y_{t-1}|I(|y_{t-1}|\leq a)\\
&\stackrel{a.s.}{\rightarrow}\phi_0\mathbb{E}(s_t-p_0)\mathbb{E}[|y_{t-1}|I(|y_{t-1}|\leq a)]
+(\mathbb{E}\varepsilon_t)\mathbb{P}(|y_{t-1}|\leq a)+0\cdot\mathbb{E}[|y_{t-1}|I(|y_{t-1}|\leq a)]=0
\end{flalign*}
and
\begin{flalign}\label{SdAsym}
\frac{1}{n}\sum_{t=1}^n(\widehat \eta_t-\bar \eta)^2&=\frac{1}{n}\sum_{t=1}^n \widehat \eta_t^2-\bar \eta^2\nonumber\\
&=\phi_0^2\frac{1}{n}\sum_{t=1}^n(s_t-p_0)^2y_{t-1}^2I(|y_{t-1}|\leq a)
+\frac{1}{n}\sum_{t=1}^n\varepsilon_t^2I(|y_{t-1}|\leq a)\nonumber\\
&\quad+(p_0\phi_0-\widehat{p}_n\widehat{\phi}_n)^2\frac{1}{n}\sum_{t=1}^ny_{t-1}^2I(|y_{t-1}|\leq a)\nonumber\\
&\quad+2\phi_0\frac{1}{n}\sum_{t=1}^n(s_t-p_0)\varepsilon_t|y_{t-1}|I(|y_{t-1}|\leq a)\nonumber\\
&\quad+2\phi_0(p_0\phi_0-\widehat{p}_n\widehat{\phi}_n)\frac{1}{n}\sum_{t=1}^n(s_t-p_0)y^2_{t-1}I(|y_{t-1}|\leq a)\\
&\quad+
2(p_0\phi_0-\widehat{p}_n\widehat{\phi}_n)\frac{1}{n}\sum_{t=1}^n\varepsilon_t|y_{t-1}|I(|y_{t-1}|\leq a)\nonumber\\
&\stackrel{a.s.}{\rightarrow}p_0(1-p_0)\phi_0^2\mathbb{E}\{y_{t-1}^2I(|y_{t-1}|\leq a)\}
+\sigma_0^2\mathbb{P}(|y_{t-1}|\leq a)=\sigma_\eta^2,\nonumber
\end{flalign}
where $\sigma_\eta^2$ is defined in (\ref{eta_sigma}). Further, using the preceding expression of $\bar{\eta}$,
we have that $\bar{\eta}=O_p(1/\sqrt{n})$ by the martingale central limit theorem
in \cite{Brown1971MartingaleCL} and Theorems \ref{thm.consistent}-\ref{thm.asym}.
Similarly, we can get
\begin{flalign*}
\frac{1}{n}\sum_{t=k+1}^n(\widehat{\eta}_{t}-\bar{\eta})(\widehat{\eta}_{t-k}-\bar{\eta})
-\frac{1}{n}\sum_{t=k+1}^n\widehat{\eta}_{t}\widehat{\eta}_{t-k}
=\frac{n-k}{n}\bar{\eta}^2
-\bar{\eta}\frac{1}{n}\sum_{t=k+1}^n\widehat{\eta}_{t}
-\bar{\eta}\frac{1}{n}\sum_{t=k+1}^n\widehat{\eta}_{t-k}
=O_p(1/n)
\end{flalign*}
for each fixed $k\geq 0$. Using above facts, we have that
\begin{flalign}\label{LI:26}
\sqrt{n}\widehat{\rho}_{nk}=(1+o_p(1))
\left\{\frac{1}{\sigma_\eta^2\sqrt{n}}\sum_{t=k+1}^n\widehat{\eta}_{t}\widehat{\eta}_{t-k}\right\}+o_p(1).
\end{flalign}

Next, it suffices to consider the joint limiting distribution of $(\sigma_\eta^2\sqrt{n})^{-1}\sum_{t=k+1}^n\widehat{\eta}_{t}\widehat{\eta}_{t-k}$, $k=1,...,M$.
To this end, let $\eta_t(\theta)=(y_t-p\phi|y_{t-1}|)I(|y_{t-1}|\leq a)$,
then $\widehat{\eta}_t=\eta_t(\widehat{\theta}_n)$ and $\eta_t=\eta_t(\theta_0)$ in (\ref{eta:expression}). Denote
\begin{flalign*}
\rho_{nk}(\theta)=\frac{1}{n\sigma_\eta^2}\sum_{t=k+1}^n\eta_t(\theta)\,\eta_{t-k}(\theta), \quad
\theta\in\Theta, \quad k\geq 1.
\end{flalign*}
Note that $\partial\eta_t(\theta)/\partial\theta=-\vartheta|y_{t-1}|I(|y_{t-1}|\leq a)$,
where $\vartheta=(p, \phi, 0)'$,
and
\begin{flalign*}
\frac{\partial\rho_{nk}(\theta_0)}{\partial\theta'}
&=\frac{1}{n\sigma_\eta^2}\sum_{t=k+1}^n\frac{\partial\eta_t(\theta_0)}{\partial\theta'}\,\eta_{t-k}
+\frac{1}{n\sigma_\eta^2}\sum_{t=k+1}^n\eta_t\frac{\partial\eta_{t-k}(\theta_0)}{\partial\theta'}\\
&=-\vartheta_0'\left\{\frac{1}{n\sigma_\eta^2}\sum_{t=k+1}^n\eta_{t-k}|y_{t-1}|I(|y_{t-1}|\leq a)
+\frac{1}{n\sigma_\eta^2}\sum_{t=k+1}^n\eta_t|y_{t-k-1}|I(|y_{t-k-1}|\leq a)\right\}\\
&=\frac{u_k}{\sigma_\eta^2}\vartheta_0'+o_p(1)
\end{flalign*}
with $u_k=-\mathbb{E}\{\eta_{t-k}|y_{t-1}|I(|y_{t-1}|\leq a)$, by the law of large numbers and $\mathbb{E}\eta_t=0$.
Then, by the Taylor expansion, the law of large numbers, and Theorems \ref{thm.consistent}-\ref{thm.asym},
it follows that
\begin{flalign*}
\sqrt{n}\big(\rho_{nk}(\widehat{\theta}_n)-\rho_{nk}(\theta_0)\big)
=\frac{\partial\rho_{nk}(\theta_0)}{\partial\theta'}\sqrt{n}(\widehat{\theta}_n-\theta_0)+o_p(1)
=\frac{u_k}{\sigma_\eta^2}\vartheta_0'\sqrt{n}(\widehat{\theta}_n-\theta_0)+o_p(1).
\end{flalign*}
Let ${\bm\rho}_n(\theta)=(\rho_{n1}(\theta), ..., \rho_{nM}(\theta))'$. It follows that
\begin{flalign*}
\sqrt{n}{\bm\rho}_n(\widehat{\theta}_n)=\sqrt{n}{\bm\rho}_n(\theta_0)+
\frac{1}{\sigma_\eta^2}\left(
                         \begin{array}{c}
                           u_1 \\
                           \vdots \\
                           u_M
                         \end{array}
                       \right)\vartheta_0'\sqrt{n}(\widehat{\theta}_n-\theta_0)+o_p(1).
\end{flalign*}
By \eqref{eq.A.4}, we have
\begin{flalign*}
\sqrt{n}(\widehat\theta_n-\theta_0)
=-\mathcal{J}^{-1}\frac{1}{\sqrt n}
\sum_{t=1}^n\frac{\partial \ell_t(\theta_0)}{\partial\theta}+o_p(1).
\end{flalign*}
The martingale central limit theorem in \cite{Brown1971MartingaleCL} gives that
\begin{flalign*}
\left(
  \begin{array}{c}
   \sqrt{n}{\bm\rho}_n(\theta_0) \\
   \sqrt{n}(\widehat\theta_n-\theta_0)
  \end{array}
\right)\stackrel{d}\longrightarrow\mathcal{N}(0, \mathbf{G}).
\end{flalign*}
Thus, $\sqrt{n}{\bm\rho}_n(\widehat{\theta}_n)\stackrel{d}\longrightarrow\
\mathcal{N}(0, \mathbf{UGU}')$ by a matrix linear transformation.

Finally, note that, by (\ref{SdAsym})-(\ref{LI:26}),
\begin{flalign*}
\sqrt{n}\big(\widehat{\bm\rho}_n-{\bm\rho}_n(\widehat{\theta}_n)\big)=o_p(1)\sqrt{n}{\bm\rho}_n(\widehat{\theta}_n)
+o_p(1)=o_p(1)O_p(1)+o_p(1)=o_p(1).
\end{flalign*}
Thus, $\sqrt{n}\widehat{\bm\rho}_n\stackrel{d}\longrightarrow\mathcal{N}(0, \mathbf{UGU}')$. The proof is
complete.\hfill$\square$

\subsection{Proof of Proposition \ref{prop:kbubblenbm}}
For any time point $t$, a $k$-th cumulative bubble collapses if $s_t = 0$, $s_{t-l} = 1$ for $1 \leq l \leq k$ and $s_{t-k-1} = 0$. Let $\{z_s^\diamond\}$ be a new auxiliary process that satisfies the recursion
\begin{flalign*}
z_s^\diamond =\left \{
\begin{array}{ll}
\varepsilon_s, &\mbox{if $s \leq t-k-1$},\\
\phi_0|z_{s-1}^\diamond| + \varepsilon_s, & \mbox{if $s > t-k-1$},
\end{array}
\right.
\end{flalign*}
then $y_{t-1} = z_{t-1}^\diamond$ if a $k$-th cumulative bubble is formed at time $t-1$ to be collapsed at time $t$. Note that the process $\{z_s^\diamond\}$ is constructed using the innovation sequence $\{\varepsilon_t\}$, which is independent of the sequence $\{s_t\}$, we can show that the joint probability
\begin{flalign*}
 & \mathbb{P}(r_t \leq c_r, s_t = 0, s_{t-1} = \cdots = s_{t-k} = 1, s_{t-k-1} = 0) \\
=\ & \mathbb{P}(\varepsilon_t - \phi_0|z_{t-1}^\diamond| \leq c_r, s_t = 0, s_{t-1} = \cdots = s_{t-k} = 1, s_{t-k-1} = 0) \\
=\ & p_0^k(1-p_0)^2 \mathbb{P}(\varepsilon_t - \phi_0|z_{t-1}^\diamond| \leq c_r).
\end{flalign*}
On the other hand, the marginal probability that a $k$-th cumulative bubble collapses at time $t$ equals to $\mathbb{P}(s_t = 0, s_{t-1} = \cdots = s_{t-k} = 1, s_{t-k-1} = 0) = p_0^k(1-p_0)^2$, and thus it suffices to show that
\begin{flalign*}
\mathbb{P}(\varepsilon_t - \phi_0|z_{t-1}^\diamond| \leq c_r) = \mathbb{P}(z_k \geq -c_r).
\end{flalign*}
For this, note that the two vectors $(\varepsilon_{t-k},\ldots,\varepsilon_t)$ and $(\varepsilon_0,\ldots,\varepsilon_k)$ share the same distribution, and thus by definition the two vectors $(z_{t-1}^\diamond,\ldots,z_{t-k-1}^\diamond)$ and $(z_k,\ldots,z_0)$ have the same joint distribution. By independence of $\varepsilon_t$ and $z_{t-1}^\diamond$ we can then conclude that
\begin{flalign*}
\mathbb{P}(\varepsilon_t - \phi_0|z_{t-1}^\diamond| \leq c_r) = \mathbb{P}(\varepsilon_k - \phi_0|z_{k-1}| \leq c_r) = \mathbb{P}(\phi_0|z_{k-1}| - \varepsilon_k \geq -c_r).
\end{flalign*}
If the innovation sequence $\{\varepsilon_t\}$ has a symmetric distribution, then $z_k = \phi_0|z_{k-1}| + \varepsilon_k$ has the same distribution as $\phi_0|z_{k-1}| - \varepsilon_k$, and the result follows. \hfill$\square$

\subsection{Proof of Proposition \ref{prop:kbubble}}
For any time point $t$, it constitutes a $k$-th cumulative bubble if $s_{t-l} = 1$ for $0 \leq l \leq k-1$ and $s_{t-k} = 0$. Let $\{z_s^\circ\}$ be a new auxiliary process that satisfies the recursion
\begin{flalign*}
z_s^\circ =\left \{
\begin{array}{ll}
\varepsilon_s, &\mbox{if $s \leq t-k$},\\
\phi_0|z_{s-1}^\circ| + \varepsilon_s, & \mbox{if $s > t-k$},
\end{array}
\right.
\end{flalign*}
then by the independence of $\{\varepsilon_t\}$ and $\{s_t\}$ we can show that the joint probability
\begin{flalign*}
\mathbb{P}(y_t > c, s_t = \cdots = s_{t-k+1} = 1, s_{t-k} = 0) & = \mathbb{P}(z_t^\circ > c, s_t = \cdots = s_{t-k+1} = 1, s_{t-k} = 0) \\
 & = \mathbb{P}(z_t^\circ > c) \cdot \mathbb{P}(s_t = \cdots = s_{t-k+1} = 1, s_{t-k} = 0) \\
 & = p_0^k(1-p_0) \mathbb{P}(z_t^\circ > c).
\end{flalign*}
On the other hand, the marginal probability that time $t$ is a $k$-th cumulative bubble equals to $\mathbb{P}(s_t = \cdots = s_{t-k+1} = 1, s_{t-k} = 0) = p_0^k(1-p_0)$, and thus it suffices to show that $\{z_s^\circ\}_{t-k < s \leq t}$ and $\{z_{s'}\}_{1 \leq s' \leq k}$ share the same distribution. For this, note that the two vectors $(\varepsilon_{t-k},\ldots,\varepsilon_t)$ and $(\varepsilon_0,\ldots,\varepsilon_k)$ share the same distribution, and they drive $z_s^\circ$, $t-k < s \leq t$, and $z_{s'}$, $1 \leq s' \leq k$, based on the same recursion, the result then follows. \hfill$\square$

%
%
%
%
%

\bibliographystyle{agsm}
\bibliography{Ref}
\end{document}